\documentclass{jams-l}

\usepackage{amssymb}
\usepackage{amsmath}

\usepackage{graphicx}

\usepackage[cmtip,all]{xy}

\usepackage[dvipsnames]{xcolor}
\usepackage{stmaryrd}
\usepackage{comment}
\usepackage{tikz}
\usepackage{caption}
\usepackage{subcaption}
\usepackage{float}
\usepackage{wrapfig}
\usepackage{bm}
\usepackage{dashbox}
\usepackage{multirow}
\usepackage{xcolor}
\usepackage{hyperref}
\usepackage{cancel}
\usepackage{mathrsfs}
\usepackage{fontawesome5}

\usepackage{todonotes}

\usetikzlibrary{quotes,angles}

\newtheorem{theorem}{Theorem}[section]
\newtheorem{lemma}[theorem]{Lemma}
\newtheorem{proposition}[theorem]{Proposition}
\newtheorem{proposition-definition}[theorem]{Proposition-Definition}
\newtheorem{corollary}[theorem]{Corollary}

\theoremstyle{definition}
\newtheorem{definition}[theorem]{Definition}

\theoremstyle{remark}
\newtheorem{remark}[theorem]{Remark}

\numberwithin{equation}{section}

\DeclareMathOperator{\vol}{Vol}

\DeclareMathOperator{\trace}{Tr}

\DeclareMathOperator{\Real}{Re}
\DeclareMathOperator{\Imag}{Im}

\DeclareMathOperator{\vect}{Vect}

\newcommand{\transpose}[1]{\phantom{}^{t} #1}

\begin{document}

\title[A Minkowski-type theorem on distances to cusps]{A Minkowski-type theorem on distances to cusps: the general case}


\author{Mathieu Dutour}
\address{Laboratoire de mathématiques de Besan\c{c}on, Université Marie et Louis Pasteur}
\curraddr{}
\email{mathieu.dutour@univ-fcomte.fr}
\thanks{}



\date{}

\dedicatory{}

\begin{abstract}
    In \cite{dutour2025minkowskitypetheoremdistancescusps}, we studied the connection between points in $\mathbb{H}^n$ and $2$-dimensional rigid adelic spaces on a totally real number field $K$ with class number $h_K = 1$. This last assumption was needed to link heights and distances to cusps. In this paper, we remove this hypothesis to obtain, without restriction on $K$ totally real, an analogue of Minkowski's second theorem on the Roy--Thunder minima of a $2$-dimensional rigid adelic space in the framework of distances between a point $\tau \in \mathbb{H}^n$ and its two closest cusps.
\end{abstract}

\maketitle

\tableofcontents


\section{Introduction}

    \subsection{Summary of the class number one case}

        The present paper is a follow-up to \cite{dutour2025minkowskitypetheoremdistancescusps}, which contains a more detailed introduction on the motivation behind the problem studied here. Let us briefly summarize the situation.

        \subsubsection{Rank $2$ lattices}

            Consider a rank $2$ Euclidean lattice $E$, which can be seen as~$\mathbb{Z}^2$ endowed with a twist of the canonical Euclidean inner-product by a symmetric, positive-definite matrix $S \in S_2^{++} \left( \mathbb{R} \right)$. Using the transformations
            \begin{equation}
            \label{eq:correspondanceS2HIntro}
                \begin{array}{ccccccc}
                    \multicolumn{3}{c}{S_2^{++} \left( \mathbb{R} \right)} & \longrightarrow & \multicolumn{3}{c}{\mathbb{H} \times \mathbb{R}_+^{\ast}} \\[1em]

                    S & = & \left( \begin{array}{cc} u & v \\[0.4em] v & w \end{array} \right) & \longmapsto & \displaystyle \left( \tau_S, \det S \right) & = & \displaystyle \left( \frac{v + i \sqrt{\det S}}{w}, \, \det S \right) \\[2em]

                    \multicolumn{3}{c}{\displaystyle \frac{\sqrt{\lambda}}{y} \left( \begin{array}{cc} x^2+y^2 & x \\[0.4em] x & 1 \end{array} \right)} & \longmapsfrom & \multicolumn{3}{c}{\displaystyle \left( \tau \, = \, x + i y, \, \lambda \right)}
                \end{array}
            \end{equation}
            we can associate to $S$ a point $\tau_S$ in the Poincaré upper half-plane $\mathbb{H}$, and this correspondence is both bijective and compatible with the natural actions of the modular group $PSL_2 \left( \mathbb{Z} \right)$. As a consequence, any property of the lattice $E$ can be seen in a geometric way on $\mathbb{H}$. Looking at the usual fundamental domain for the action of $PSL_2 \left( \mathbb{Z} \right)$ on $\mathbb{H}$ by homographic transformations

            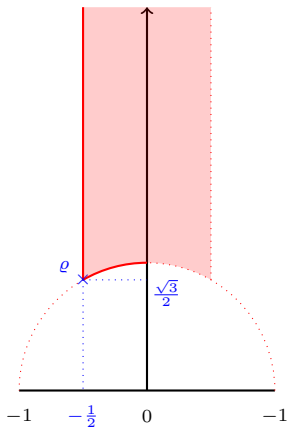
\begin{figure}[H]
    
            \centering
        
                \begin{tikzpicture}[scale=1.7]
        
                    \draw[thick] (-1,0) coordinate (1) node {} -- (1,0) coordinate (2) node {};
                    \draw[thick, ->] (0,0) -- (0,3);
                
                    \node (A) at (-0.5, 0.8660254) {};
                    \node (B) at (0.5,0.8660254) {};
                
                    \draw[thick,red] (0,1) arc (90:120:1cm);
                    \draw[dotted,red] (1,0) arc (0:90:1cm);
                
                    \draw[thick,red] (-0.5,0.8660254) -- (-0.5,3);
                    \draw[dotted,red] (0.5,0.8660254) -- (0.5,3);
                
                    \draw[dotted,red] (-0.5,0.8660254) arc (120:180:1cm);

                    \draw[dotted,blue] (-1/2,0.866025) -- (0,0.866025);
                    \draw[dotted,blue] (-1/2,0.866025) -- (-1/2,0);
                
                    \node (C) at (-1,-0.2) {$\scriptstyle -1$};
                    \node (D) at (1,-0.2) {$\scriptstyle -1$};
                    \node (E) at (0,-0.2) {$\scriptstyle 0$};
                    \node (E) at (-1/2,0.866025) {$\scriptstyle \color{blue} \times$};
                    \node (F) at (-1/2-0.15,0.866025+0.1) {$\scriptstyle \color{blue} \varrho$};
                    \node (G) at (0.15,0.866025-0.1) {$\scriptstyle \color{blue} \frac{\sqrt{3}}{2}$};
                    \node (H) at (-1/2,-0.2) {$\scriptstyle \color{blue} -\frac{1}{2}$};
                
                    \fill[red,opacity=0.2] (A) arc (120:60:1cm) -- (0.5,3) -- (-0.5,3) -- (A) -- cycle;

                \end{tikzpicture}
        
            \caption{Fundamental domain for the action of $PSL_2 \left( \mathbb{Z} \right)$ on $\mathbb{H}$}
            \label{fig:fundamentalDomainPSL2Z}
    
            \end{figure}

            \noindent one sees that the point which is the farthest away from the cusp $\infty$ corresponds to the rank $2$ lattice whose first minimum (\textit{i.e.} the shortest length of a non-zero vector) is maximal.

        \subsubsection{Rigid adelic spaces in the class number one case}

            In \cite{dutour2025minkowskitypetheoremdistancescusps}, we generalized the observation made above. Consider a totally real number field $K$ of degree $n$. By iterating \eqref{eq:correspondanceS2HIntro}, one can attach a $2$-dimensional rigid adelic space $E_{\tau}$ with total height $1$ to any point $\tau \in \mathbb{H}^n$. One can define, for a $2$-dimensional rigid adelic space~$E$, the \textit{Roy--Thunder minima} of $E$ as
            \begin{equation}
                \begin{array}{lll}
                    \Lambda_1 \left( E \right) & = & \inf \left \{ H_{E} \left( x \right) \; \middle \vert \; x \in E \setminus \left \{ 0 \right \} \right \}, \\[0.5em]

                    \Lambda_2 \left( E \right) & = & \inf \left \{ \max \left( H_{E} \left( x \right), H_{} \left( y \right) \right) \; \middle \vert \; \vect_K \left( x,y \right) \, = \, E \right \},
                \end{array}
            \end{equation}
            and a Hermite-type constant
            \begin{equation}
            \label{eq:intro:HermiteConstant}
                \begin{array}{lll}
                    \displaystyle c_{II}^{\Lambda} \left( 2, K \right) & = & \displaystyle \sup\limits_{\dim E = 2} \sqrt{\frac{\Lambda_1 \left( E \right) \Lambda_2 \left( E \right)}{H \left( E \right)}}
                \end{array}
            \end{equation}
            where $H \left( E \right)$ is the height of $E$. The correspondence between points in $\mathbb{H}^n$ and rigid adelic spaces allowed us, assuming that the class number $h_K$ of $K$ is $1$, to interpret~$\Lambda_1 \left( E_{\tau} \right)$ and $\Lambda_2 \left( E_{\tau} \right)$ as powers of the distances between $\tau$ and its two closest cusps, and to obtain the following version of Minkowski's second theorem.
            
            \begin{theorem}
                For any $\tau \in \mathbb{H}^n$, we have
                \begin{equation}
                    \begin{array}{lllll}
                        \displaystyle \frac{1}{c_{II}^{\Lambda} \left( 2, K \right)^{4n}} & \leqslant & \mu_1 \left( \tau \right) \mu_2 \left( \tau \right) & \leqslant & 1,
                    \end{array}
                \end{equation}
                where $c_{II}^{\Lambda} \left( 2, K \right)$ is the constant \eqref{eq:intro:HermiteConstant}.
            \end{theorem}
            \noindent Let us insist that the assumption $h_K=1$ was needed there. As a consequence, we were able, in this situation, to make improvements on \cite[Lemmas I.2.1, I.2.2]{vanDerGeer:hilbert-modular-surfaces}.

    \subsection{Statement of the results}

        This paper is devoted to generalizing \cite{dutour2025minkowskitypetheoremdistancescusps} to account for all totally real number fields, without restriction on their class number. Consider such a number field $K$ of degree $n$, whose ring of integers is denoted by~$\mathcal{O}_K$. Let $\mathfrak{a}$ be a fractional ideal of $\mathcal{O}_K$. To any point $\tau \in \mathbb{H}^n$, we attach in definition \ref{def:rigidAdelicSpaceEATau} a $2$-dimensional rigid adelic space $E_{\mathfrak{a}, \tau}$ for which the \textit{Roy--Thunder minima} equal
        \begin{equation}
            \begin{array}{lll}
                \Lambda_1 \left( E_{\mathfrak{a}, \tau} \right) & = & N \left( \mathfrak{a} \right)^{-1/n} \mu_{\mathfrak{a}, 1} \left( \tau \right)^{-1/2n}, \\[0.5em]

                \Lambda_2 \left( E_{\mathfrak{a}, \tau} \right) & = & N \left( \mathfrak{a} \right)^{-1/n} \mu_{\mathfrak{a}, 2} \left( \tau \right)^{-1/2n},
            \end{array}
        \end{equation}
        where the functions $\mu_{\mathfrak{a},1} \left( \tau \right)^{-1/2}$ and $\mu_{\mathfrak{a},2} \left( \tau \right)^{-1/2}$ represent the $\mathfrak{a}$-distances between~$\tau$ and its two closest cusps. The analogue of Minkowski's second theorem on these distances takes the following form.
        \begin{theorem}
            For any fractional ideal $\mathfrak{a}$ and any $\tau \in \mathbb{H}^n$, we have
            \begin{equation}
                \begin{array}{lllll}
                    \displaystyle \frac{1}{c_{II}^{\Lambda} \left( 2, K \right)^{4n}} \cdot \frac{1}{N \left( \mathfrak{a} \right)^2} & \leqslant & \mu_{\mathfrak{a}, 1} \left( \tau \right) \mu_{\mathfrak{a}, 2} \left( \tau \right) & \leqslant & \displaystyle \frac{1}{N \left( \mathfrak{a} \right)^2},
                \end{array}
            \end{equation}
            where $c_{II}^{\Lambda} \left( 2, K \right)$ is the constant introduced in \cite[Definition 2.31]{dutour2025minkowskitypetheoremdistancescusps} and in \eqref{eq:intro:HermiteConstant}.
        \end{theorem}

        This theorem has several consequences. The first two, which generalize similar observations made in \cite{dutour2025minkowskitypetheoremdistancescusps}, are improvements on results in \cite{vanDerGeer:hilbert-modular-surfaces}, this time without restriction on the class number of $K$.
        \begin{enumerate}
            \item Corollary \ref{cor:separationCusps} is an effective and uniform version of \cite[Lemma I.2.1]{vanDerGeer:hilbert-modular-surfaces};
            \item Corollary \ref{cor:lowerBoundMu1} is an optimal version of \cite[Lemma I.2.2]{vanDerGeer:hilbert-modular-surfaces}
        \end{enumerate}
        Once again, using the proof of \cite[Lemma I.2.2]{vanDerGeer:hilbert-modular-surfaces}, and the optimal nature of corollary~\ref{cor:lowerBoundMu1} with $\mathfrak{a} = \mathcal{O}_K$, we obtain
        \begin{equation}
            \begin{array}{lll}
                c_{II}^{\Lambda} \left( 2, K \right) & \leqslant & \sqrt{2} \Delta_K^{1/2n},
            \end{array}
        \end{equation}
        \textit{i.e.} the upper-bound obtained by Gaudron and Rémond in \cite[Proposition 5.1]{gaudron-remond:corps-siegel}. In the last section of this paper, we give several statements, whose proofs are entirely similar to the ones present in \cite[Section 4.5]{dutour2025minkowskitypetheoremdistancescusps}, culminating in the following result.

        \begin{theorem}
                For any $0 \leqslant t < 1$, we have
                \begin{equation}
                    \begin{array}{lllll}
                        \multicolumn{5}{l}{\displaystyle \frac{1}{N \left( \mathfrak{a} \right)^t} \left( \frac{1}{c_{II}^{\Lambda} \left( 2, K \right)^{2nt}} \left( 1 - \frac{1}{c_{II}^{\Lambda} \left( 2, K \right)^{2n}} \right) + \frac{1}{1-t} \cdot \frac{1}{c_{II}^{\Lambda} \left( 2, K \right)^{2n}} \right)} \\[2em]

                        & \leqslant & \displaystyle \frac{1}{\vol \left( \widehat{\Gamma}_K \left( \mathfrak{a} \right) \backslash \mathbb{H}^n \right)} \int_{\widehat{\Gamma}_K \left( \mathfrak{a} \right) \backslash \mathbb{H}^n} \mu_{\mathfrak{a},1} \left( \tau \right)^t \; \mathrm{d}m\left( \tau \right) & \leqslant & \displaystyle \frac{1}{1-t} \cdot \frac{1}{N \left( \mathfrak{a} \right)^t},
                    \end{array}
                \end{equation}
                where $\mathrm{d}m \left( \tau \right)$ denotes the Poincaré measure on $\mathbb{H}^n$.
            \end{theorem}

    \subsection{Future work}

        In the introduction of \cite{dutour2025minkowskitypetheoremdistancescusps}, the section on future work contained several questions, the first of them being answered in this paper. Let us provide updates on the other problems posed therein.

        \begin{enumerate}
            \item Several geometric properties seen in this paper could be useful in obtaining more information on the Hermite-type constants $c_{II}^{\Lambda} \left( 2, K \right)$.
                \begin{enumerate}
                    \item On the sphere of influence of a cusp $c$, the point $\tau_c$ which is the farthest away from $c$ is at $\mathfrak{a}$-distance
                    \begin{equation}
                        \begin{array}{lll}
                            \mu_{\mathfrak{a},1} \left( \tau_c \right)^{-1/2} & = & c_{II}^{\Lambda} \left( 2, K \right)^n N \left( \mathfrak{a} \right)^{1/2}
                        \end{array}
                    \end{equation}
                    from $c$. Whether $\tau_c$ has any further interesting properties, such as being elliptical, is not known.

                    \item It seems possible that, on a given sphere of influence, only finitely many points are at a given $\mathfrak{a}$-distance from their two closest cusps. If one was able to interpret these distances as coordinates, perhaps a computation of volume would yield a result in terms of the Hermite-type constant which, using the classical formula for the volume of a Hilbert modular variety, would provide a closed formula for $c_{II}^{\Lambda} \left( 2, K \right)$. Whether this can be done is unknown.
                \end{enumerate}

            \item Dealing with rigid adelic spaces, or even lattices, of greater rank poses significant challenges. The most natural analogues of $\mathbb{H}$ and $PSL_2 \left( \mathbb{Z} \right)$ in this setting would be the Siegel upper half-space $\mathcal{H}_g$ and the Siegel modular group $Sp_{2g} \left( \mathbb{Z} \right)$. However, the analogue of identification \eqref{eq:correspondanceS2HIntro} is only a bijection
            \begin{equation}
            \label{eq:higherRankCorrespondence}
                \begin{array}{ccc}
                    \mathcal{H}_g & \overset{\sim}{\longrightarrow} & S_{2g}^{++} \left( \mathbb{R} \right) \cap Sp_{2g} \left( \mathbb{R} \right),
                \end{array}
            \end{equation}
            which is compatible with the action of $Sp_{2g} \left( \mathbb{Z} \right)$, but does not extend with this compatibility to a map whose image would be $S_{2g}^{++} \left( \mathbb{R} \right)$. Whether we can use \eqref{eq:higherRankCorrespondence} to obtain interesting results related to what is presented in this paper remains to be seen.

            \item Finally, having the analogue of Minkowski's theorem on the distances to cusps without restriction on the totally real number field $K$ opens the door to an application in the situation presented in \cite{frey-lefourn-lorenzo:height-estimates, habegger-pazuki:bad-reduction}. In particular, one could aim to obtain an answer to \cite[Question 4.1]{frey-lefourn-lorenzo:height-estimates}, though it seems unlikely that a uniform lower-bound could be obtained using the tools presented in this paper.
        \end{enumerate}

    \subsection{Acknowledgements}

        I would like to thank Samuel Le Fourn for discussions around this topic, and for providing invaluable background on the questions which motivated this project.

\section{Algebraic number theory}

    \subsection{The ideal class groups}

        Compared to \cite[Section 2]{dutour2025minkowskitypetheoremdistancescusps}, we need to go in more details regarding the ideal class group and its narrow version. Let us consider a totally real number field $K$, whose ring of integers is denoted by $\mathcal{O}_K$.

        \begin{definition}
            A \textit{fractional ideal} $\mathfrak{a}$ of $\mathcal{O}_K$ is a non-zero, finitely generated $\mathcal{O}_K$-submodule of $K$. The set of fractional ideals of $\mathcal{O}_K$ is denoted by $J_K$.
        \end{definition}

        \begin{proposition}
            The set $J_K$ is an Abelian group for the product of ideals, with identity element $\mathcal{O}_K$. The inverse of any $\mathfrak{a} \in J_K$ is further given by
            \begin{equation}
                \begin{array}{lll}
                    \mathfrak{a}^{-1} & = & \left \{ x \in K \; \middle \vert \; x \mathfrak{a} \subseteq \mathcal{O}_K \right \}.
                \end{array}
            \end{equation}
        \end{proposition}

        \begin{proof}
            This is given in \cite[Proposition I.3.8, p.21]{neukirch:algebraic-number-theory}.
        \end{proof}

        \begin{definition}
            A fractional ideal $\mathfrak{a}$ is said to be \textit{principal} if it is of the form $x \mathcal{O}_K$ with $x \in K^{\ast}$. Their set is a subgroup of $J_K$, denoted by $P_K$.
        \end{definition}

        Let us denote by $\sigma_1, \ldots, \sigma_n$ the embeddings of $K$ in $\mathbb{R}$.

        \begin{definition}
            An element $x \in K^{\ast}$ is called \textit{totally positive} if we have
            \begin{equation}
                \begin{array}{llll}
                    \sigma_j \left( x \right) & > & 0 & \text{for any } j \in \llbracket 1, n \rrbracket.
                \end{array}
            \end{equation}
            The set of totally positive elements of $K$ forms a subgroup of $K^{\ast}$, which we denote by $K^+$.
        \end{definition}

        \begin{definition}
            A fractional ideal $\mathfrak{a}$ is said to be \textit{totally positive principal} if it is of the form $x \mathcal{O}_K$ with $x \in K^{+}$. Their set $P_K^+$ is a subgroup of~$J_K$.
        \end{definition}

        \begin{definition}
            The \textit{ideal class group} is defined as the quotient
            \begin{equation}
                \begin{array}{lll}
                    Cl_K & = & J_K / P_K
                \end{array}
            \end{equation}
            and the \textit{narrow ideal class group} as the quotient
            \begin{equation}
                \begin{array}{lll}
                    Cl_K^+ & = & J_K / P_K^+.
                \end{array}
            \end{equation}
        \end{definition}

        \begin{remark}
            There is a natural embedding $Cl_K \hookrightarrow Cl_K^+$, coming from the inclusion of $P_K^+$ in $P_K$.
        \end{remark}

        \begin{proposition}
            The ideal class group $Cl_K$ is finite. Its order is denoted by $h_K$, and called the class number of $K$.
        \end{proposition}

        \begin{proof}
            The finiteness of the ideal class group is proved in \cite[Theorem I.6.3, p.36]{neukirch:algebraic-number-theory}. It actually requires the notion of absolute norms of ideals, which will be recalled in the next subsection.
        \end{proof}

        \begin{proposition}
            There is an exact sequence
            \begin{equation}
                \begin{array}{ccccccccccc}
                    1 & \longrightarrow & \mathcal{O}_K^{\times}/\mathcal{O}_K^{\times,+} & \longrightarrow & K^{\ast}/K^+ & \longrightarrow & Cl_K^+ & \longrightarrow & Cl_K & \longrightarrow & 1
                \end{array}
            \end{equation}
            of groups. Consequently, the narrow ideal class group $Cl_K^+$ is finite, and its order~$h_K^+$ is given by
            \begin{equation}
                \begin{array}{lll}
                    h_K^+ & = & \displaystyle h_K \cdot \frac{[ K^{\ast} : K^+ ]}{[ \mathcal{O}_K^{\times} : \mathcal{O}_K^{\times,+}]},
                \end{array}
            \end{equation}
            where the quotient of indices on the right-hand side is a power of two.
        \end{proposition}

        \begin{proof}
            This is a direct computation, noting that the map $Cl_K^+ \longrightarrow Cl_K$ is the projection map, and that the map
            \begin{equation}
                \begin{array}{ccc}
                    K^{\ast} & \longrightarrow & Cl_K^+ \\[0.5em]
                    x & \longmapsto & x \mathcal{O}_K
                \end{array}
            \end{equation}
            induces the map $K^{\ast}/K^{\times} \longrightarrow Cl_K^+$.
        \end{proof}

        \begin{remark}
            The image of a fractional ideal $\mathfrak{a}$ in $Cl_K^+$ is called its \textit{narrow class}.
        \end{remark}

        \begin{remark}
            Since $\mathcal{O}_K$ is a \textit{Dedekind domain}, every fractional ideal $\mathfrak{a}$ is generated by (at most) two elements. Namely, we can find $\alpha, \beta \in K$ such that we have
            \begin{equation}
                \begin{array}{lll}
                    \mathfrak{a} & = & \alpha \mathcal{O}_K + \beta \mathcal{O}_K.
                \end{array}
            \end{equation}
            There is a stronger version of this result, which is that if $\mathfrak{a}, \mathfrak{b}$ are two fractional ideals, there exists $\alpha, \beta \in K$ such that we have
            \begin{equation}
                \begin{array}{lll}
                    \mathfrak{a} & = & \alpha \mathcal{O}_K + \beta \mathfrak{b}.
                \end{array}
            \end{equation}
        \end{remark}

        We conclude this subsection by defining a particular fractional ideal of $\mathcal{O}_K$ which will be used later.

        \begin{definition}
        \label{def:differentIdeal}
            The \textit{codifferent ideal} is defined as
            \begin{equation}
                \begin{array}{lll}
                    \mathfrak{c}_K & = & \left \{ x \in K \; \middle \vert \; \forall y \in \mathcal{O}_K, \; \trace \left( xy \right) \in \mathbb{Z} \right \},
                \end{array}
            \end{equation}
            the trace of an $u$ element of $K$ being defined by
            \begin{equation}
                \begin{array}{lll}
                    \trace \left( u \right) & = & \sum\limits_{j=1}^n \sigma_j \left( u \right).
                \end{array}
            \end{equation}
            The \textit{different ideal} $\mathfrak{d}_K$ is then defined as
            \begin{equation}
                \begin{array}{lll}
                    \mathfrak{d}_K & = & \mathfrak{c}_K^{-1}.
                \end{array}
            \end{equation}
        \end{definition}

        \begin{remark}
            The fractional ideal $\mathfrak{c}_K$ contains~$\mathcal{O}_K$, and its inverse $\mathfrak{d}_K$ is contained in $\mathcal{O}_K$. Furthermore, though we have $\mathfrak{d}_K^{-1} = \mathfrak{c}_K$, the notation $\mathfrak{d}_K^{-1}$ is still often used.
        \end{remark}

    \subsection{Norm of an ideal}

        \begin{definition}
            The (absolute) \textit{ideal norm} of an ideal $\mathfrak{a}$ of $\mathcal{O}_K$ is defined as
            \begin{equation}
                \begin{array}{lll}
                    N \left( \mathfrak{a} \right) & = & \# \mathcal{O}_K / \mathfrak{a}.
                \end{array}
            \end{equation}
        \end{definition}

        \begin{proposition}
            For any two ideals $\mathfrak{a}, \mathfrak{b}$ of $\mathcal{O}_K$, we have
            \begin{equation}
                \begin{array}{lll}
                    N \left( \mathfrak{a} \mathfrak{b} \right) & = & N \left( \mathfrak{a} \right) N \left( \mathfrak{b} \right).
                \end{array}
            \end{equation}
            Consequently, the norm map extends to a map
            \begin{equation}
                \begin{array}{ccccc}
                    N & : & J_K & \longrightarrow & \mathbb{Q}_+^{\ast},
                \end{array}
            \end{equation}
            where we recall that $J_K$ is the group of fractional ideals.
        \end{proposition}

        \begin{proof}
            The multiplicativity of the absolute norm of ideals is a direct consequence of the Chinese remainder theorem, and of the fact that (fractional) ideals of $\mathcal{O}_K$ may be written (uniquely) as product of prime ideals. This is explained in more details in \cite[Proposition I.6.1]{neukirch:algebraic-number-theory}.
        \end{proof}

        \begin{proposition}
            The norm of the different ideal $\mathfrak{d}_K$ is given by
            \begin{equation}
                \begin{array}{lll}
                    N \left( \mathfrak{d}_K \right) & = & \Delta_K.
                \end{array}
            \end{equation}
        \end{proposition}

        \begin{proof}
            This is a particular case of \cite[Theorem 2.9, p.201]{neukirch:algebraic-number-theory}.
        \end{proof}

        \noindent The aim of this paragraph is to provide a formula computing $N \left( \mathfrak{a} \right)$ in terms of the generators of $\mathfrak{a}$ as an ideal of $\mathcal{O}_K$.

        \begin{proposition}
        \label{prop:isoProdLocalizations}
            There exists a ring isomorphism
            \begin{equation}
                \begin{array}{lll}
                    \mathcal{O}_K / \mathfrak{a} & \simeq & \prod\limits_{v \in V_f \left( K \right)} \left( \mathcal{O}_K \right)_{\mathfrak{p}_v} / \mathfrak{a} \left( \mathcal{O}_K \right)_{\mathfrak{p}_v}.
                \end{array}
            \end{equation}
        \end{proposition}

        \begin{proof}
            The following argument is a detailed version of \cite[Proposition I.12.3]{neukirch:algebraic-number-theory}. We begin by setting, for any $v \in V_f \left( K \right)$,
            \begin{equation}
                \begin{array}{lll}
                    \mathfrak{a}_v & = & \mathcal{O}_K \cap \mathfrak{a} \left( \mathcal{O}_K \right)_{\mathfrak{p_v}}.
                \end{array}
            \end{equation}
            Since $\mathcal{O}_K$ is a Dedekind domain, the ideal $\mathfrak{a}$ can be written in a unique way as a product of (finitely many of) the prime ideals $\mathfrak{p}_v$. As a consequence, we have $\mathfrak{a} \not \subset \mathfrak{p}_v$ for all but finitely many places $v \in V_f \left( K \right)$. For such places, we have
            \begin{equation}
                \begin{array}{lll}
                    \mathfrak{a} \left( \mathcal{O}_K \right)_{\mathfrak{p}_v} & = & \left( \mathcal{O}_K \right)_{\mathfrak{p}_v},
                \end{array}
            \end{equation}
            since elements which are in $\mathfrak{a}$ but not in $\mathfrak{p}_v$ are invertible in $\left( \mathcal{O}_K \right)_{\mathfrak{p}_v}$, and thus we have $\mathfrak{a}_v = \mathcal{O}_K$. Therefore, we have
            \begin{equation}
                \begin{array}[t]{lll}
                    \bigcap\limits_{v \in V_f \left( K \right)} \mathfrak{a}_v & = & \bigcap\limits_{\substack{v \in V_f \left( K \right) \\ \mathfrak{a} \subset \mathfrak{p}_v}} \mathfrak{a}_v,
                \end{array}
            \end{equation}
            and we will now prove that this intersection equals $\mathfrak{a}$. Consider an element
            \begin{equation}
                \begin{array}[t]{lllll}
                    a & \in & \bigcap\limits_{v \in V_f \left( K \right)} \mathfrak{a}_v & = & \bigcap\limits_{\substack{v \in V_f \left( K \right) \\ \mathfrak{a} \subset \mathfrak{p}_v}} \mathfrak{a}_v,
                \end{array}
            \end{equation}
            and set
            \begin{equation}
                \begin{array}{lll}
                    \mathfrak{b} & = & \left \{ x \in \mathcal{O}_K \; \middle \vert \; xa \in \mathfrak{a} \right \},
                \end{array}
            \end{equation}
            which is an ideal of $\mathcal{O}_K$. Let $v \in V_f \left( K \right)$. We have
            \begin{equation}
                \begin{array}{lll}
                    a & \in & \mathfrak{a} \left( \mathcal{O}_K \right)_{\mathfrak{p}_v},
                \end{array}
            \end{equation}
            so there exists $s \in \mathcal{O}_K \setminus \mathfrak{p}_v$ such that we have $sa \in \mathfrak{a}$. This gives $s \in \mathfrak{b}$, which yields~$\mathfrak{b} \not \subset \mathfrak{p}_v$. The ideals $\mathfrak{p}_v$ being the only maximal ideals of $\mathcal{O}_K$, we get
            \begin{equation}
                \begin{array}{lll}
                    \mathfrak{b} & = & \mathcal{O}_K.
                \end{array}
            \end{equation}
            We thus have $1 \in \mathfrak{b}$, which is to say we have $a \in \mathfrak{a}$, thereby giving the inclusion
            \begin{equation}
            \label{eq:inclusionIntersectionInA}
                \begin{array}[t]{lll}
                    \bigcap\limits_{\substack{v \in V_f \left( K \right) \\ \mathfrak{a} \subset \mathfrak{p}_v}} \mathfrak{a}_v & \subseteq & \mathfrak{a}.
                \end{array}
            \end{equation}
            The converse inclusion to \eqref{eq:inclusionIntersectionInA} is direct. Once more, let $v \in V_f \left( K \right)$, and assume we have $\mathfrak{a} \subseteq \mathfrak{p}_v$. The radical of the ideal $\mathfrak{a}_v$, defined as
            \begin{equation}
                \begin{array}{lll}
                    \sqrt{\mathfrak{a}_v} & = & \left \{ x \in \mathcal{O}_K \; \middle \vert \; \exists n \in \mathbb{N}, \; x^n \in \mathfrak{a}_v \right \},
                \end{array}
            \end{equation}
            is then a proper ideal of $\mathcal{O}_K$, or we would have $1 \in \mathcal{O}_K \cap \mathfrak{a} \left( \mathcal{O}_K \right)_{\mathfrak{p}_v}$, and there would exist $s \in \mathcal{O}_K \setminus\mathfrak{p}_v$ with $s \cdot 1 \in \mathfrak{a}$, contradicting the assumption~$\mathfrak{a} \subseteq \mathfrak{p}_v$. The localized ring $\left( \mathcal{O}_K \right)_{\mathfrak{p}_v}$ being of discrete valuation, there exists $n \in \mathbb{N}^{\ast}$ such that we have
            \begin{equation}
                \begin{array}{lll}
                    \mathfrak{a} \left( \mathcal{O}_K \right)_{\mathfrak{p}_v} & = & \mathfrak{p}_v^n \left( \mathcal{O}_K \right)_{\mathfrak{p}_v}.
                \end{array}
            \end{equation}
            We then have
            \begin{equation}
                \begin{array}{lllllll}
                    \mathfrak{p}_v & = & \sqrt{\mathfrak{p}_v^n} & = & \sqrt{\mathfrak{p}_v^n \left( \mathcal{O}_K \right)_{\mathfrak{p}_v} \cap \mathcal{O}_K} & \subseteq & \sqrt{\mathfrak{a}_v}.
                \end{array}
            \end{equation}
            Since $\sqrt{\mathfrak{a}_v}$ is a proper ideal of $\mathcal{O}_K$, we get\footnote{Though unnecessary in this paper, the ideals $\mathfrak{a}_v$ provide a \textit{primary decomposition} of $\mathfrak{a}$.}
            \begin{equation}
                \begin{array}{lll}
                    \sqrt{\mathfrak{a}_v} & = & \mathfrak{p}_v.
                \end{array}
            \end{equation}
            As a consequence, the only maximal ideal of $\mathcal{O}_K$ containing $\mathfrak{a}_v$ is precisely $\mathfrak{p}_v$. To see this, consider $v' \in V_f \left( K \right)$ a place of $K$ satisfying $\mathfrak{a}_{v} \subseteq \mathfrak{p}_{v'}$. We have
            \begin{equation}
                \begin{array}{lllllll}
                    \mathfrak{p}_v & = & \sqrt{\mathfrak{a}_v} & \subseteq & \sqrt{\mathfrak{p}_{v'}} & = & \mathfrak{p}_{v'},
                \end{array}
            \end{equation}
            which yields $v=v'$. Thus, if $v, v' \in V_f \left( K \right)$ are distinct places of $K$, we must have
            \begin{equation}
                \begin{array}{lll}
                    \mathfrak{a}_v + \mathfrak{a}_{v'} & = & \mathcal{O}_K,
                \end{array}
            \end{equation}
            so the Chinese remainder theorem (see \cite[Section I.3.6]{neukirch:algebraic-number-theory}) applies, and yields a ring isomorphism
            \begin{equation}
            \label{eq:CRTAv}
                \begin{array}{lll}
                    \mathcal{O}_K / \mathfrak{a} & \simeq & \prod\limits_{v \in V_f \left( K \right)} \mathcal{O}_K / \mathfrak{a}_v.
                \end{array}
            \end{equation}
            Consider a place $v \in V_f \left( K \right)$. Using again the fact that $\mathfrak{p}_v$ is the only maximal ideal of $\mathcal{O}_K$ containing $\mathfrak{a}_v$, we see that the quotient ring $\mathcal{O}_K / \mathfrak{a}_v$ is local, thereby giving
            \begin{equation}
            \label{eq:quotientByAv}
                \begin{array}{llllll}
                    \mathcal{O}_K / \mathfrak{a}_v & = & \left( \mathcal{O}_K / \mathfrak{a}_v \right)_{\mathfrak{p}_v} & \simeq & \left( \mathcal{O}_K \right)_{\mathfrak{p}_v} / \mathfrak{a}_v \left( \mathcal{O}_K \right)_{\mathfrak{p}_v},
                \end{array}
            \end{equation}
            since taking the quotient commutes with localizing. The equality
            \begin{equation}
                \begin{array}{lll}
                    \mathfrak{a}_v \left( \mathcal{O}_K \right)_{\mathfrak{p}_v} & = & \mathfrak{a} \left( \mathcal{O}_K \right)_{\mathfrak{p}_v}
                \end{array}
            \end{equation}
            can be plugged into \eqref{eq:quotientByAv}, giving a ring isomorphism
            \begin{equation}
            \label{eq:isoQuotientByAv}
                \begin{array}{lll}
                    \mathcal{O}_K / \mathfrak{a}_v & \simeq & \left( \mathcal{O}_K \right)_{\mathfrak{p}_v} / \mathfrak{a} \left( \mathcal{O}_K \right)_{\mathfrak{p}_v}.
                \end{array}
            \end{equation}
            Combining the ring isomorphisms \eqref{eq:CRTAv} and \eqref{eq:isoQuotientByAv} concludes the proof.
        \end{proof}

        \begin{corollary}
        \label{cor:isoQuotientOkACompletions}
            There exists a ring isomorphism
            \begin{equation}
                \begin{array}{lll}
                    \mathcal{O}_K / \mathfrak{a} & \simeq & \prod\limits_{v \in V_f \left( K \right)} \mathcal{O}_{K_v} / \mathfrak{a} \mathcal{O}_{K_v}.
                \end{array}
            \end{equation}
        \end{corollary}

        \begin{proof}
            Using proposition \ref{prop:isoProdLocalizations}, we need only prove that, for any place $v \in V_f \left( K \right)$, there exists a ring isomorphism
            \begin{equation}
                \begin{array}{lll}
                    \left( \mathcal{O}_K \right)_{\mathfrak{p}_v} / \mathfrak{a} \left( \mathcal{O}_K \right)_{\mathfrak{p}_v} & \simeq & \mathcal{O}_{K_v} / \mathfrak{a} \mathcal{O}_{K_v}.
                \end{array}
            \end{equation}
            To that effect, let us adapt the argument presented in \cite[Proposition II.2.4]{neukirch:algebraic-number-theory}, as hinted at by Neukirch in \cite[Proposition II.4.3]{neukirch:algebraic-number-theory}. Consider a place $v \in V_f \left( K \right)$. Let us recall that the ring of integers
            \begin{equation}
                \begin{array}{lll}
                    \mathcal{O}_{K_v} & = & \left \{ x \in K_v \; \middle \vert \; \left \vert x \right \vert_v \, \leqslant \, 1 \right \}
                \end{array}
            \end{equation}
            is the completion of $\left( \mathcal{O}_K \right)_{\mathfrak{p}_v}$ with respect to $\left \vert \cdot \right \vert_v$. The ring $\left( \mathcal{O}_K \right)_{\mathfrak{p}_v}$ being of discrete valuation, there exists $n \in \mathbb{N}$ such that we have
            \begin{equation}
                \begin{array}{lll}
                    \mathfrak{a} \left( \mathcal{O}_K \right)_{\mathfrak{p}_v} & = & \mathfrak{p}_v^n \left( \mathcal{O}_K \right)_{\mathfrak{p}_v}.
                \end{array}
            \end{equation}
            By completion, we also have
            \begin{equation}
                \begin{array}{lll}
                    \mathfrak{a} \mathcal{O}_{K_v} & = & \mathfrak{p}_v^n \mathcal{O}_{K_v}.
                \end{array}
            \end{equation}
            Embedding $\left( \mathcal{O}_K \right)_{\mathfrak{p}_v}$ into its completion $\mathcal{O}_{K_v}$ and projecting modulo $\mathfrak{p}_v^n \mathcal{O}_{K_v}$, we get a ring homomorphism
            \begin{equation}
            \label{eq:mapEmbedProjectCompletion}
                \begin{array}{lll}
                    \left( \mathcal{O}_K \right)_{\mathfrak{p}_v} & \longrightarrow & \mathcal{O}_{K_v} / \mathfrak{p}_v^n \mathcal{O}_{K_v},
                \end{array}
            \end{equation}
            which we will now prove is surjective. Let $x \in \mathcal{O}_{K_v}$. By completion, there exists a sequence of elements $u_n \in \left( \mathcal{O}_K \right)_{\mathfrak{p}_v}$ such that we have 
            \begin{equation}
                \begin{array}{lll}
                    v \left( x - u_n \right) & \geqslant & n.
                \end{array}
            \end{equation}
            for any $n \in \mathbb{N}$. In other words, we have
            \begin{equation}
                \begin{array}{lll}
                    x - u_n & \in & \mathfrak{p}_v^n \mathcal{O}_{K_v},
                \end{array}
            \end{equation}
            so the map \eqref{eq:mapEmbedProjectCompletion} is indeed surjective. Let us study its kernel. Consider $x \in \left( \mathcal{O}_K \right)_{\mathfrak{p}_v}$ satisfying $x \in \mathfrak{p}_v^n \mathcal{O}_{K_v}$. We have $v \left( x \right) \geqslant n$, which gives $x \in \mathfrak{p}_v^n \left( \mathcal{O}_K \right)_{\mathfrak{p}_v}$. Thus \eqref{eq:mapEmbedProjectCompletion} induces a ring isomorphism
            \begin{equation}
                \begin{array}{lllllll}
                    \left( \mathcal{O}_K \right)_{\mathfrak{p}_v} / \mathfrak{a} \left( \mathcal{O}_K \right)_{\mathfrak{p}_v} & = & \left( \mathcal{O}_K \right)_{\mathfrak{p}_v} / \mathfrak{p}_v^n \left( \mathcal{O}_K \right)_{\mathfrak{p}_v} \\[1em]
                    
                    & \simeq & \mathcal{O}_{K_v} / \mathfrak{p}_v^n \mathcal{O}_{K_v} & = & \mathcal{O}_{K_v} / \mathfrak{a} \mathcal{O}_{K_v},
                \end{array}
            \end{equation}
            and completes the proof of the proposition.
        \end{proof}

        \begin{proposition}
        \label{prop:numberElementsOkvXv}
            For any $v \in V_f \left( K \right)$, and any non-zero $x_v \in \mathcal{O}_{K_v}$, we have
            \begin{equation}
                \begin{array}{lll}
                    \# \mathcal{O}_{K_v}/x_v \mathcal{O}_{K_v} & = & \left \vert x_v \right \vert_v^{-n_v}.
                \end{array}
            \end{equation}
        \end{proposition}

        \begin{proof}
            Let $v \in V_f \left( K \right)$ and $x_v \in \mathcal{O}_{K_v}$. Denote by $p$ the prime number such that we have $v \, \vert \, p$. Recall that $\mathcal{O}_{K_v}$ is a free $\mathbb{Z}_p$-module of rank $n_v = \left[ K_v : \mathbb{Q}_p \right]$. Since~$\mathbb{Z}_p$ is a principal ideal domain, we can find a $\mathbb{Z}_p$-basis $\left( e_1, \ldots, e_{n_v} \right)$ of $\mathcal{O}_{K_v}$ and $p$-adic integers $d_1$, \ldots, $d_{n_v} \in \mathbb{Z}_p \setminus \left \{ 0 \right \}$ such that $\left( d_1 e_1, \ldots, d_{n_v} e_{n_v} \right)$ is a $\mathbb{Z}_p$-basis of $x_v \mathcal{O}_{K_v}$. We then have
            \begin{equation}
                \begin{array}{lll}
                    \mathcal{O}_{K_v} / x_v \mathcal{O}_{K_v} & \simeq & \prod\limits_{j=1}^{n_v} \mathbb{Z}_p / d_j \mathbb{Z}_p.
                \end{array}
            \end{equation}
            Without loss of generality, we may choose $d_j$ to be given by $d_j = p^{m_j}$ with $m_j \geqslant 0$, thus giving
            \begin{equation}
                \begin{array}{lllll}
                    \# \mathbb{Z}_p / p^{m_j} \mathbb{Z}_p & = & p^{m_j} & = & \left \vert d_j \right \vert_p^{-1}.
                \end{array}
            \end{equation}
            Consequently, we have
            \begin{equation}
                \begin{array}{lll}
                    \# \mathcal{O}_{K_v} / x_v \mathcal{O}_{K_v} & = & \prod\limits_{j=1}^{n_v} \left \vert d_j \right \vert_p^{-1}.
                \end{array}
            \end{equation}
            Let us now consider the $\mathbb{Q}_p$-linear map
            \begin{equation}
                \begin{array}{ccccc}
                    \varphi & : & K_v & \longrightarrow & K_v \\
                    && e_j & \longmapsto & d_j e_j
                \end{array}.
            \end{equation}
            Its determinant is given by
            \begin{equation}
                \begin{array}{lll}
                    \det \varphi & = & \prod\limits_{j=1}^{n_v} d_j,
                \end{array}
            \end{equation}
            and it sends $\mathcal{O}_{K_v}$ onto $x_v \mathcal{O}_{K_v}$. Denote by $\psi : K_v \longrightarrow K_v$ the $\mathbb{Q}_p$-linear map induced by the multiplication by $x_v$. By definition, we have
            \begin{equation}
                \begin{array}{lll}
                    \det \psi & = & N_{K_v/\mathbb{Q}_p} \left( x_v \right).
                \end{array}
            \end{equation}
            These maps $\varphi$ and $\psi$ are invertible, and we can set
            \begin{equation}
                \begin{array}{lllllll}
                    \chi & = & \varphi \circ \psi^{-1} & : & K_v & \longrightarrow & K_v.
                \end{array}
            \end{equation}
            Since $\chi$ induces an automorphism of the $\mathbb{Z}_p$-module $x_v \mathcal{O}_{K_v}$, we have
            \begin{equation}
                \begin{array}{lllll}
                    \frac{\det \varphi}{\det \psi} & = & \det \chi & \in & \mathbb{Z}_p^{\times},
                \end{array}
            \end{equation}
            which gives
            \begin{equation}
            \label{eq:prodDjNormXv}
                \begin{array}{lllllll}
                    \prod\limits_{j=1}^{n_v} d_j & = & \det \varphi & = & \varepsilon \det \psi & = & \varepsilon N_{K_v/\mathbb{Q}_p} \left( x_v \right),
                \end{array}
            \end{equation}
            with $\varepsilon \in \mathbb{Z}_p^{\times}$. Applying the $p$-adic norms to both sides of \eqref{eq:prodDjNormXv} and then taking the inverses yield
            \begin{equation}
                \begin{array}{lllllll}
                    \# \mathcal{O}_{K_v} / x_v \mathcal{O}_{K_v} & = & \prod\limits_{j=1}^{n_v} \left \vert d_j \right \vert_p^{-1} & = & \left \vert N_{K_v/\mathbb{Q}_p} \left( x_v \right) \right \vert_p^{-1} & = & \left \vert x_v \right \vert_v^{-n_v},
                \end{array}
            \end{equation}
            thus concluding the proof.
        \end{proof}

        \begin{corollary}
        \label{cor:computeNormIdeal}
            Let $\alpha, \beta \in K$ and $\mathfrak{a}$ be a fractional ideal of $\mathcal{O}_K$. Assume the fractional ideal $\mathfrak{b} = \alpha \mathcal{O}_K + \beta \mathfrak{a}$ of $\mathcal{O}_K$ is non-zero. For any $v \in V_f \left( K \right)$, set $x_v \in K_v$ such that we have $\mathfrak{a} \mathcal{O}_{K_v} \, = \, x_v \mathcal{O}_{K_v}$. We have
            \begin{equation}
            \label{eq:normOfIdeal}
                \begin{array}[t]{lll}
                    N \left( \mathfrak{b} \right) & = & \prod\limits_{v \in V_f \left( K \right)} \max \left( \left \vert \alpha \right \vert_v, \left \vert \beta x_v \right \vert_v \right)^{-n_v}.
                \end{array}
            \end{equation}
        \end{corollary}

        \begin{proof}
            Let us first consider the case where $\mathfrak{a} \subseteq \mathcal{O}_K$ is an integral ideal, and $\alpha, \beta$ are in $\mathcal{O}_K$. Consider a place $v \in V_f \left( K \right)$. The ring $\mathcal{O}_{K_v}$ being a principal ideal domain, there does indeed exist elements $x_v, y_v \in \mathcal{O}_{K_v}$ such that we have
            \begin{equation}
                \begin{array}{lllllll}
                    & \mathfrak{a} \mathcal{O}_{K_v} & = & x_v \mathcal{O}_{K_v} \\[0.5em]
                    \text{and} & \mathfrak{b} \mathcal{O}_{K_v} & = & y_v \mathcal{O}_{K_v} & = & \alpha \mathcal{O}_{K_v} + \beta x_v \mathcal{O}_{K_v}.
                \end{array}
            \end{equation}
            These elements $x_v$ and $y_v$ further satisfy
            \begin{equation}
            \label{eq:maxNormsAlphaBeta}
                \begin{array}{lll}
                    \left \vert y_v \right \vert_v & = & \max \left( \left \vert \alpha \right \vert_v, \left \vert \beta x_v \right \vert_v \right).
                \end{array}
            \end{equation}
            We have
            \begin{equation}
            \label{eq:normIdealB}
                \begin{array}{llll}
                    N \left( \mathfrak{b} \right) & = & \prod\limits_{v \in V_f \left( K \right)} \# \mathcal{O}_{K_v} / \mathfrak{b} \mathcal{O}_{K_v} & \text{by corollary \ref{cor:isoQuotientOkACompletions}} \\[2em]

                    & = & \prod\limits_{v \in V_f \left( K \right)} \# \mathcal{O}_{K_v} / y_v \mathcal{O}_{K_v} \\[2em]

                    & = & \prod\limits_{v \in V_f \left( K \right)} \left \vert y_v \right \vert_v^{-n_v} & \text{by proposition \ref{prop:numberElementsOkvXv}} \\[2em]

                    & = & \prod\limits_{v \in V_f \left( K \right)} \max \left( \left \vert \alpha \right \vert_v, \left \vert \beta x_v \right \vert_v \right)^{-n_v} & \text{by \eqref{eq:maxNormsAlphaBeta}.}
                \end{array}
            \end{equation}
            The case $\alpha, \beta \in K$ is dealt with by clearing out denominators and using the product formula. Finally, let us consider the situation in which $\mathfrak{a}$ is a fractional ideal. We set $u \in \mathcal{O}_K$ non-zero such that we have $u \mathfrak{a} \subseteq \mathcal{O}_K$. We have
            \begin{equation}
                \begin{array}{llll}
                    N \left( \mathfrak{b} \right) & = & N_{K/\mathbb{Q}} \left( u \right)^{-1} N \left( \alpha u \mathcal{O}_K + \beta \left( u \mathfrak{a} \right) \right) \\[1em]

                    & = & N_{K/\mathbb{Q}} \left( u \right)^{-1} \prod\limits_{v \in V_f \left( K \right)} \max \left( \left \vert \alpha u \right \vert_v, \left \vert \beta u x_v \right \vert_v \right)^{-n_v} & \text{by \eqref{eq:normIdealB}} \\[2em]

                    & = & \prod\limits_{v \in V_f \left( K \right)} \max \left( \left \vert \alpha \right \vert_v, \left \vert \beta x_v \right \vert_v \right)^{-n_v},
                \end{array}
            \end{equation}
            which concludes the proof.
        \end{proof}

        \begin{remark}
            Taking $\mathfrak{a} \, = \, \mathcal{O}_K$ in corollary \ref{cor:computeNormIdeal}, and using the fact that any fractional ideal $\mathfrak{b}$ of~$\mathcal{O}_K$ can be generated by (at most) two elements $\alpha, \beta \in K$, we get, if $\mathfrak{b}$ is non-zero,
            \begin{equation}
                \begin{array}{lll}
                    N \left( \mathfrak{b} \right) & = & \prod\limits_{v \in V_f \left( K \right)} \max \left( \left \vert \alpha \right \vert_v, \left \vert \beta \right \vert_v \right)^{-n_v}.
                \end{array}
            \end{equation}
        \end{remark}

\section{\texorpdfstring{Action of generalized Hilbert modular groups on $\mathbb{H}^n$}{Action of generalized Hilbert modular groups on Hn}}

    Much of the material we will need regarding Hilbert modular groups and varieties has already been presented in \cite[Section 3]{dutour2025minkowskitypetheoremdistancescusps}. Thus, we will only add a few elements to this previous overview. Throughout this section, we will consider by $K$ a totally real number field of degree $n$, whose ring of integers is denoted by $\mathcal{O}_K$, and whose real embeddings are $\sigma_1$, \ldots, $\sigma_n$.

    \subsection{Generalized Hilbert modular groups}

        Let $\mathfrak{a}$ be a fractional ideal of $\mathcal{O}_K$.

        \begin{definition}
            The subring $M_{2, \mathfrak{a}}$ of $M_2 \left( K \right)$ is defined as
            \begin{equation}
                \begin{array}{lll}
                    M_{2,\mathfrak{a}} & = & \left \{ \left( \begin{array}{cc} a & b \\[0.4em] c & d \end{array} \right) \; \middle \vert \; a,d \in \mathcal{O}_K, \; b \in \mathfrak{a}^{-1}, \; c \in \mathfrak{a} \right \}.
                \end{array}
            \end{equation}
        \end{definition}

        \begin{definition}
        \label{def:generalizedHMG}
            The \textit{generalized Hilbert modular group} $\widehat{\Gamma}_K \left( \mathfrak{a} \right)$ is defined as
            \begin{equation}
                \begin{array}{lll}
                    \widehat{\Gamma}_K \left( \mathfrak{a} \right) & = & \left \{ M \in M_{2,\mathfrak{a}} \; \middle \vert \; \det M \in \mathcal{O}_K^{\times,+} \right\} / \left \{ \left( \begin{array}{cc} \varepsilon & 0 \\[0.4em] 0 & \varepsilon \end{array} \right) \; \middle \vert \; \varepsilon \in \mathcal{O}_K^{\times} \right \}.
                \end{array}
            \end{equation}
        \end{definition}

        \begin{remark}
            Note that taking $\mathfrak{a} \, = \, \mathcal{O}_K$ in \ref{def:generalizedHMG} gives the Hilbert modular group $\widehat{\Gamma}_K$ presented in \cite[Section 3.1.2]{dutour2025minkowskitypetheoremdistancescusps}. Using the real embeddings $\sigma_1$, \ldots, $\sigma_n$, we further note that $\widehat{\Gamma}_K \left( \mathfrak{a} \right)$ is embedded in $PGL_2 \left( \mathbb{R} \right)^n$, and thus acts on $\mathbb{H}^n$.
        \end{remark}

        \begin{proposition}
        \label{prop:multByLambdaMap}
            Consider a totally positive element $\lambda \in K$. The map
            \begin{equation}
            \label{eq:mapHMVNarrowClass}
                \begin{array}{ccccc}
                    \varphi & : & \widehat{\Gamma}_K \left( \lambda \mathfrak{a} \right) \backslash \mathbb{H}^n & \longrightarrow & \widehat{\Gamma}_K \left( \mathfrak{a} \right) \backslash \mathbb{H}^n \\[1em]
                    && \tau & \longmapsto & \lambda \tau
                \end{array}
            \end{equation}
            is a biholomorphism which induces the equality
            \begin{equation}
                \begin{array}{lll}
                    \vol \left( \widehat{\Gamma}_K \left( \lambda \mathfrak{a} \right) \backslash \mathbb{H}^n \right) & = & \vol \left( \widehat{\Gamma}_K \left( \mathfrak{a} \right) \backslash \mathbb{H}^n \right)
                \end{array}
            \end{equation}
            with respect to the Poincaré metric. In other words, the volume of the Hilbert modular variety $\widehat{\Gamma}_K \left( \mathfrak{a} \right) \backslash \mathbb{H}^n$ only depends on the narrow class of $\mathfrak{a}$.
        \end{proposition}

        \begin{proof}
            Let us begin by noting that we set
            \begin{equation}
                \begin{array}{lllll}
                    \lambda \tau & = & \lambda \left( \tau_1, \ldots, \tau_n \right) & = & \left( \sigma_1 \left( \lambda \right) \tau_1, \ldots, \sigma_n \left( \lambda \right) \tau_n \right).
                \end{array}
            \end{equation}
            Since $\lambda$ is a totally positive element of $K$, multiplication by $\lambda$ is a biholomorphism of $\mathbb{H}^n$, and \eqref{eq:mapHMVNarrowClass} is a well-defined holomorphic map, as it satisfies
            \begin{equation}
                \begin{array}{lll}
                    \varphi \left( \left( \begin{array}{cc} a & \lambda^{-1}b \\[0.4em] \lambda c & d \end{array} \right) \cdot \tau \right) & = & \lambda \left( \begin{array}{cc} a & b \\[0.4em] c & d \end{array} \right) \cdot \tau, 
                \end{array}
            \end{equation}
            It is furthermore bijective, and its holomorphic inverse is given by multiplication by $1/\lambda$. The equality of volume can be seen using the explicit description of the measure
            \begin{equation}
                \begin{array}{lll}
                \text{d}m \left( \tau \right) & = & \displaystyle \frac{\text{d}x_1 \text{d}y_1 \ldots \text{d}x_n \text{d}y_n}{y_1^2 \cdots y_n^2}
                \end{array}
            \end{equation}
            of the measure associated with the Poincaré metric.
        \end{proof}

    \subsection{Cusps}

        As in \cite[Section 2.1]{dutour2025minkowskitypetheoremdistancescusps}, the action of $\widehat{\Gamma}_K \left( \mathfrak{a} \right)$ leads to the definition of the cusps. In what follows, we consider a fractional ideal $\mathfrak{a}$ of $\mathcal{O}_K$.

        \subsubsection{The set of cusps}

            Even though this paragraph mirrors \cite[Section 2.1.1]{dutour2025minkowskitypetheoremdistancescusps} to a great extent, it will be useful to make notations precise for the rest of this text.

            \begin{proposition}
                The set of cusps for the action of $\widehat{\Gamma}_K \left( \mathfrak{a} \right)$, defined as
                \begin{equation}
                    \begin{array}{lll}
                        \left \{ c \in \mathbb{P}^1 \left( \mathbb{R} \right)^n \; \middle \vert \; \exists \gamma \in \widehat{\Gamma}_K \left( \mathfrak{a} \right), \; \gamma \cdot c \, = \, c \right \}
                    \end{array}
                \end{equation}
                is equal to (and identified with) the embedding of $\mathbb{P}^1 \left( K \right)$ in $\mathbb{P}^1 \left( \mathbb{R} \right)^n$ using the real embeddings $\sigma_1$, \ldots, $\sigma_n$. 
            \end{proposition}

            \begin{proof}
                This is a direct computation.
            \end{proof}

            \begin{definition}
                For any cusp $c \in \mathbb{P}^1 \left( K \right)$, we define
                \begin{equation}
                    \begin{array}{lll}
                        \widehat{\Gamma}_K \left( \mathfrak{a} \right)_c & = & \left \{ \gamma \in \widehat{\Gamma}_K \left( \mathfrak{a} \right) \; \middle \vert \; \gamma \cdot c \, = \, c \right \}
                    \end{array}
                \end{equation}
                the stabilizer of $c$ in $\widehat{\Gamma}_K \left( \mathfrak{a} \right)$.
            \end{definition}

            \begin{proposition}
                The well-defined map
                \begin{equation}
                    \begin{array}{ccccc}
                        \widehat{\Gamma}_K \left( \mathfrak{a} \right) \backslash \mathbb{P}^1 \left( K \right) & \longrightarrow & Cl_K \\[0.5em]
                        \left[ \alpha : \beta \right] & \longmapsto & \alpha \mathcal{O}_K + \beta \mathfrak{a}^{-1}
                    \end{array}
                \end{equation}
                is bijective. Thus, there are exactly $h_K$ cusps up to $\widehat{\Gamma}_K \left( \mathfrak{a} \right)$-equivalence.
            \end{proposition}

            \begin{proof}
                The proof is similar to the one used in \cite[Proposition I.1.1]{vanDerGeer:hilbert-modular-surfaces}, with some minor modifications to account for the presence of the fractional ideal $\mathfrak{a}$.
            \end{proof}

            \begin{definition}
                The cusp $\infty$ is defined as $\left[ 1 : 0 \right] \in \mathbb{P}^1 \left( K \right)$.
            \end{definition}

            \begin{proposition}
                For any fractional ideal $\mathfrak{b}$ of $\mathcal{O}_K$, the stabilizer $\widehat{\Gamma}_K \left( \mathfrak{b} \right)_{\infty}$ of the cusp $\infty$ is given by
                \begin{equation}
                    \begin{array}{lll}
                        \widehat{\Gamma}_K \left( \mathfrak{b} \right)_{\infty} & = & \left \{ \left( \begin{array}{cc} \varepsilon & \mu \\[0.4em] 0 & 1 \end{array} \right) \; \middle \vert \; \varepsilon \in \mathcal{O}_K^{\times, +}, \; \mu \in \mathfrak{b}^{-1} \right \}.
                    \end{array}
                \end{equation}
            \end{proposition}

            \begin{proof}
                This is done in \cite[Section I.4]{vanDerGeer:hilbert-modular-surfaces}.
            \end{proof}

            \begin{proposition}
            \label{prop:conjugStabilizerCusp}
                Let $c = \left[ \alpha : \beta \right] \in \mathbb{P}^1 \left( K \right)$ be a cusp, with $\alpha, \beta \in \mathcal{O}_K$. There exists $M \in SL_2 \left( K \right)$ such that we have $M \cdot \infty = c$, and
                \begin{equation}
                    \begin{array}{lll}
                        \widehat{\Gamma}_K \left( \mathfrak{a} \right) & = & M \widehat{\Gamma}_K \left( \mathfrak{a} \mathfrak{q}^2 \right) M^{-1},
                    \end{array}
                \end{equation}
                where the fractional ideal $\mathfrak{q}$ is defined by
                \begin{equation}
                    \begin{array}{lll}
                        \mathfrak{q} & = & \alpha \mathcal{O}_K + \beta \mathfrak{a}^{-1}.
                    \end{array}
                \end{equation}
                In particular, we have
                \begin{equation}
                    \begin{array}{lll}
                        \widehat{\Gamma}_K \left( \mathfrak{a} \right)_c & = & M \widehat{\Gamma}_K \left( \mathfrak{a} \mathfrak{q}^2 \right)_{\infty} M^{-1}.
                    \end{array}
                \end{equation}
            \end{proposition}

            \begin{proof}
                Let us begin by noting that we have
                \begin{equation}
                    \begin{array}{lllll}
                        \mathcal{O}_K & = & \mathfrak{q} \mathfrak{q}^{-1} & = & \alpha \mathfrak{q}^{-1} + \beta \mathfrak{a}^{-1} \mathfrak{q}^{-1}.
                    \end{array}
                \end{equation}
                Hence, we can choose $\beta^{\ast} \in \mathfrak{q}^{-1}$ and $\alpha^{\ast} \in \mathfrak{a}^{-1} \mathfrak{q}^{-1}$ such that we have
                \begin{equation}
                    \begin{array}{lll}
                        \alpha \beta^{\ast} - \beta \alpha^{\ast} & = & 1.
                    \end{array}
                \end{equation}
                Let us now consider the matrix
                \begin{equation}
                    \begin{array}{lllll}
                        M & = & \left( \begin{array}{cc} \alpha & \alpha^{\ast} \\[0.4em] \beta & \beta^{\ast} \end{array} \right) & \in & SL_2 \left( K \right),
                    \end{array}
                \end{equation}
                which satisfies
                \begin{equation}
                    \begin{array}{lllll}
                        M \cdot \infty & = & \left[ \alpha : \beta \right] & = & c.
                    \end{array}
                \end{equation}
                For any $a,d \in \mathcal{O}_K$ and any $b \in \mathfrak{a}^{-1}$, $c \in \mathfrak{a}$, we have
                \begin{equation}
                    \begin{array}{llll}
                        \multicolumn{2}{l}{M^{-1} \left( \begin{array}{cc} a & b \\[0.4em] c & d \end{array} \right) M} & = & \left( \begin{array}{cc} \beta^{\ast} & - \alpha^{\ast} \\[0.4em] -\beta & \alpha \end{array} \right) \left( \begin{array}{cc} a & b \\[0.4em] c & d \end{array} \right) \left( \begin{array}{cc} \alpha & \alpha^{\ast} \\[0.4em] \beta & \beta^{\ast} \end{array} \right) \\[2em]

                        = \hspace{-6pt} & \multicolumn{3}{l}{\left( \begin{array}{cc} a \alpha \beta^{\ast} - c \alpha \alpha^{\ast} + b \beta \beta^{\ast} - d \beta \alpha^{\ast} & \left( a - d \right) \alpha^{\ast} \beta^{\ast} - c \left( \alpha^{\ast} \right)^2 + b \left( \beta^{\ast} \right)^2 \\[0.4em] c \alpha^2 + \left( d - \alpha \right) \alpha \beta - d \beta^{2} & c \alpha \alpha^{\ast} - a \beta \alpha^{\ast} + d \alpha \beta^{\ast} - b \beta \beta^{\ast} \end{array} \right),}
                    \end{array}
                \end{equation}
                and this matrix is seen to belong to $\widehat{\Gamma}_K \left( \mathfrak{a} \mathfrak{q}^2 \right)$.
            \end{proof}

            The remainder of this paragraph will be devoted to finding a fundamental domain for the action of $\widehat{\Gamma}_K \left( \mathfrak{b} \right)_{\infty}$ on $\mathbb{H}^n$. Using proposition \ref{prop:conjugStabilizerCusp}, we will then also have found such a domain for the action of the stabilizer $\widehat{\Gamma}_K \left( \mathfrak{a} \right)_c$. As in \cite{dutour2025minkowskitypetheoremdistancescusps}, we will set
            \begin{equation}
                \begin{array}{ccccc}
                    \Real \tau & = & \left( \Real \tau_1, \ldots, \Real \tau_n \right) & \in & \mathbb{R}^n, \\[0.5em]

                    \Imag \tau & = & \left( \Imag \tau_1, \ldots, \Imag \tau_n \right) & \in & \left( \mathbb{R}_+^{\ast} \right)^n,
                \end{array}
            \end{equation}
            to lighten notations.

            \begin{proposition}
            \label{prop:fundamentalDomainStabInfty}
                A fundamental domain for the action of $\widehat{\Gamma}_K \left( \mathfrak{b} \right)_{\infty}$ on $\mathbb{H}^n$ is given by the set
                \begin{equation}
                    \begin{array}{lllll}
                        E & = & \left \{ \tau \in \mathbb{H}^n \; \middle \vert \; \Real \tau \in T, \; \Imag \tau \in F \right \} & = & T \times F,
                    \end{array}
                \end{equation}
                where $T$ is a fundamental domain for the action\footnote{Similar to \cite[Proposition 2.11]{dutour2025minkowskitypetheoremdistancescusps}, using the embeddings $\sigma_1$, \ldots, $\sigma_n$.} of $\mathfrak{b}^{-1}$ on $\mathbb{R}^n$, and $F$ is a fundamental domain for the action\footnote{See \cite[Equation $\left( 2.26 \right)$]{dutour2025minkowskitypetheoremdistancescusps}} of $\mathcal{O}_K^{\times,+}$ on $\left(\mathbb{R}_+^{\ast} \right)^n$.
            \end{proposition}

            \begin{proof}
                The proof follows the same argument as \cite[Proposition 3.17]{dutour2025minkowskitypetheoremdistancescusps}.
            \end{proof}

            An example of fundamental domain $F$ we can plug into proposition \ref{prop:fundamentalDomainStabInfty} can be found in \cite[Proposition 3.19]{dutour2025minkowskitypetheoremdistancescusps}. The appropriate results are stated here for clarity.

            \begin{lemma}
                The map
                \begin{equation}
                    \begin{array}{ccccc}
                        \psi & : & \left( \mathbb{R}_+^{\ast} \right)^n & \longrightarrow & \left( \mathbb{R}_+^{\ast} \right)^n \\[0.5em]

                        && \left( y_1, \ldots, y_n \right) & \longmapsto & \left( y_1, \ldots, y_{n-1}, y_1 \cdot \ldots \cdot y_n \right)
                    \end{array}
                \end{equation}
                is a diffeomorphism which satisfies $\psi \left( B \right) \, = \, \left( \mathbb{R}_+^{\ast} \right)^{n-1} \times \left\{ 1 \right \}$, and whose Jacobian at $\left( y_1, \ldots, y_n \right)$ equals $y_1 \cdot \ldots \cdot y_{n-1}$. Furthermore, we have
                \begin{equation}
                    \begin{array}{lll}
                        \psi \left( \varepsilon \cdot \left( y_1, \ldots, y_n \right) \right) & = & \left( \varepsilon \cdot \left( y_1, \ldots, y_{n-1} \right), y_1 \cdot \ldots \cdot y_n \right).
                    \end{array}
                \end{equation}
            \end{lemma}

            \begin{proof}
                This is a direct computation.
            \end{proof}

            \begin{proposition}
                Denote by $G$ a fundamental domain\footnote{See \cite[Proposition 2.19]{dutour2025minkowskitypetheoremdistancescusps}} for the lattice $\mathcal{O}_K^{\times,+}$ embedded in $\mathbb{R}^{n-1}$. The set
                \begin{equation}
                \label{eq:fundamentalDomainF}
                    \begin{array}{lll}
                        F & = & \psi^{-1} \left( \lambda^{-1} \left( G \right) \times \mathbb{R}_+^{\ast} \right)
                    \end{array}
                \end{equation}
                is a fundamental domain for the action of $\mathcal{O}_K^{\times,+}$ on $\left( \mathbb{R}_+^{\ast} \right)^n$. For any $r>0$, we further have
                \begin{equation}
                \label{eq:fundamentalDomainFr}
                    \begin{array}{lll}
                        F_r & = & \psi^{-1} \left( \lambda^{-1} \left( G \right) \times \left] 1/r^2, + \infty \right[ \right) \\[1em]

                        & = & \left \{ \left( y_1, \ldots, y_n \right) \in \left( \mathbb{R}_+^{\ast} \right)^n \; \middle \vert \; \left( y_1 \cdot \ldots \cdot y_n \right)^{-1/2} \, < \, r \right \}.
                    \end{array}
                \end{equation}
            \end{proposition}

            \begin{proof}
                This is a direct computation.
            \end{proof}

        \subsubsection{\texorpdfstring{The $\mathfrak{a}$-distance to the cusps}{The a-distance to the cusps}}

            In order to measure the distance between a point in $\mathbb{H}^n$ and a cusp, the traditional $\mu$-function presented in \cite{dutour2025minkowskitypetheoremdistancescusps, vanDerGeer:hilbert-modular-surfaces} does not work, as we now allow the presence of a fractional ideal in the generalized Hilbert modular groups. Consider a fractional ideal $\mathfrak{a}$ of $\mathcal{O}_K$.

            \begin{definition}
            \label{def:aDistance}
                Let $c = \left[ \alpha : \beta \right] \in \mathbb{P}^1 \left( K \right)$ be a cusp, with $\alpha, \beta \in \mathcal{O}_K$. We define
                \begin{equation}
                    \begin{array}{ccccc}
                        \mu_{\mathfrak{a}} \left( \cdot, c \right) & : & \mathbb{H}^n & \longrightarrow & \mathbb{R}_+^{\ast} \\[1em]
                        && \tau & \longmapsto & \displaystyle \frac{N \left( \alpha \mathcal{O}_K + \beta \mathfrak{a}^{-1} \right)^2 N \left( \Imag \tau \right)}{\left \vert N \left( - \beta \tau + \alpha \right) \right \vert^2}
                    \end{array},
                \end{equation}
                with $N \left( w \right) = w_1 \ldots w_n$ for any $w = \left( w_1, \ldots, w_n \right) \in \mathbb{C}^n$. The function
                \begin{equation}
                    \begin{array}{ccccc}
                        \mu_{\mathfrak{a}} \left( \cdot, c \right)^{-1/2} & : & \mathbb{H}^n & \longrightarrow &\mathbb{R}_+^{\ast}
                    \end{array}
                \end{equation}
                then represents the $\mathfrak{a}$\textit{-distance to the cusp $c$}.
            \end{definition}

            \begin{remark}
            \label{rmk:cuspInftyClassicalMu}
                Taking $c = \infty$ yields $\mu_{\mathfrak{a}} = \mu$, \textit{i.e.} the classical $\mu$-function presented in \cite[Section 3.2.2]{dutour2025minkowskitypetheoremdistancescusps}, regardless of the fractional ideal $\mathfrak{a}$.
            \end{remark}

            As evidenced by the proposition below, the function $\mu_{\mathfrak{a}}$ is compatible with the action of $\widehat{\Gamma}_K \left( \mathfrak{a} \right)$.

            \begin{proposition}
            \label{prop:functionMuActionGamma}
                For any $\gamma \in \widehat{\Gamma}_K \left( \mathfrak{a} \right)$, any $c \in \mathbb{P}^1 \left( K \right)$, and any $\tau \in \mathbb{H}^n$, we have
                \begin{equation}
                    \begin{array}{lll}
                        \mu_{\mathfrak{a}} \left( \gamma \cdot \tau, \gamma \cdot c \right) & = & \mu_{\mathfrak{a}} \left( \tau, c \right).
                    \end{array}
                \end{equation}
            \end{proposition}

            \begin{proof}
                This is a direct computation.
            \end{proof}

            Another compatibility property is to understand how $\mu_{\mathfrak{a}}$ depends on the representative of the narrow ideal class of $\mathfrak{a}$.

            \begin{proposition}
            \label{prop:muFunctionNarrowClass}
                Let $\lambda \in K$ be a totally positive element. For any $\tau \in \mathbb{H}^n$ and any cusp $c = \left[ \alpha : \beta \right] \in \mathbb{P}^1 \left( K \right)$, we have
                \begin{equation}
                    \begin{array}{lll}
                        \mu_{\lambda \mathfrak{a}} \left( \tau, \left[ \alpha : \beta \right] \right) & = & \displaystyle \frac{1}{N_{K/\mathbb{Q}} \left( \lambda \right)} \mu_{\mathfrak{a}} \left( \lambda \tau, \left[ \alpha : \lambda^{-1} \beta \right] \right).
                    \end{array}
                \end{equation}
            \end{proposition}

            \begin{proof}
                Let $\tau$ and $c$ as in the statement of this proposition. We have
                \begin{equation}
                    \begin{array}{lllll}
                        \multicolumn{3}{l}{\mu_{\lambda \mathfrak{a}} \left( \tau, \left[ \alpha : \beta \right] \right)} & = & \displaystyle N \left( \alpha \mathcal{O}_K + \beta \lambda^{-1} \mathfrak{a}^{-1} \right) \prod\limits_{j=1}^n \frac{y_j}{\left \vert - \sigma_j \left( \beta \right) \tau_j + \sigma_j \left( \alpha \right) \right \vert^2} \\[2em]

                        & = & \multicolumn{3}{l}{\displaystyle \frac{1}{N_{K/\mathbb{Q} \left( \lambda \right)}} N \left( \alpha \mathcal{O}_K + \beta \lambda^{-1} \mathfrak{a}^{-1} \right) \prod\limits_{j=1}^n \frac{\sigma_j \left( \lambda \right) y_j}{\left \vert - \sigma_j \left( \lambda^{-1} \beta \right)\left( \lambda \tau \right)_j + \sigma_j \left( \alpha \right) \right \vert^2}} \\[2em]

                        & = & \multicolumn{3}{l}{\displaystyle \frac{1}{N_{K/\mathbb{Q}} \left( \lambda \right)} \mu_{\mathfrak{a}} \left( \lambda \tau, \left[ \alpha : \lambda^{-1} \beta \right] \right).}
                    \end{array}
                \end{equation}
            \end{proof}

            Let us see one last of these properties.

            \begin{proposition}
            \label{prop:distanceMatrixMCToInfitnity}
                Consider a cusp $c = \left[ \alpha : \beta \right] \in \mathbb{P}^1 \left( K \right)$, and, as in proposition~\ref{prop:conjugStabilizerCusp}, a matrix $M \in SL_2 \left( K \right)$ such that we have $M \cdot \infty = c$ and
                \begin{equation}
                    \begin{array}{lll}
                        \widehat{\Gamma}_K \left( \mathfrak{a} \right) & = & M \widehat{\Gamma}_K \left( \mathfrak{a} \mathfrak{q}^2 \right) M^{-1},
                    \end{array}
                \end{equation}
                where we have set
                \begin{equation}
                    \begin{array}{lll}
                        \mathfrak{q} & = & \alpha \mathcal{O}_K + \beta \mathfrak{a}^{-1}.
                    \end{array}
                \end{equation}
                For any point $\tau \in \mathbb{H}^n$, we have
                \begin{equation}
                    \begin{array}{lll}
                        \mu_{\mathfrak{a}} \left( M \cdot \tau, c \right) & = & N \left( \mathfrak{q} \right)^2 \mu \left( \tau, \infty \right).
                    \end{array}
                \end{equation}
            \end{proposition}

            \begin{proof}
                Let us write the matrix $M$ as
                \begin{equation}
                    \begin{array}{lll}
                        M & = & \left( \begin{array}{cc} \alpha & \alpha^{\ast} \\[0.4em] \beta & \beta^{\ast} \end{array} \right).
                    \end{array}
                \end{equation}
                For any point $\tau \in \mathbb{H}^n$, we have
                \begin{equation}
                    \begin{array}{lll}
                        N \left( \Imag \left( M \cdot \tau \right) \right) & = & \displaystyle \frac{N \left( \Imag \tau \right)}{\left \vert N \left( \beta \tau + \beta^{\ast} \right) \right \vert^2},
                    \end{array}
                \end{equation}
                and we also have
                \begin{equation}
                    \begin{array}{lllll}
                        - \beta M \cdot \tau + \alpha & = & \displaystyle - \frac{\alpha \beta \tau + \alpha^{\ast} \beta}{\beta \tau + \beta^{\ast}} + \alpha & = & \displaystyle \frac{1}{\beta \tau + \beta^{\ast}}.
                    \end{array}
                \end{equation}
                Thus, we have
                \begin{equation}
                    \begin{array}{lll}
                        \mu_{\mathfrak{a}} \left( M \cdot \tau, c \right) & = &  \displaystyle \frac{N \left( \alpha \mathcal{O}_K + \beta \mathfrak{a}^{-1} \right)^2 N \left( \Imag \left( M \cdot \tau \right) \right)}{\left \vert N \left( - \beta M \cdot \tau + \alpha \right) \right \vert^2} \\[2em]

                        & = & N \left( \mathfrak{q} \right)^2 N \left( \Imag \tau \right) \\[1em]

                        & = & N \left( \mathfrak{q} \right)^2 \mu \left( \tau, \infty \right).
                    \end{array}
                \end{equation}
            \end{proof}

            \begin{definition}
                For any $\tau \in \mathbb{H}^n$, we set
                \begin{equation}
                \label{eq:defMu1}
                    \begin{array}{lll}
                        \mu_{\mathfrak{a},1} \left( \tau \right) & = & \max\limits_{c \in \mathbb{P}^1 \left( K \right)} \mu_{\mathfrak{a}} \left( \tau, \, c \right)
                    \end{array}
                \end{equation}
                and, denoting by $c_{\tau} \in \mathbb{P}^1 \left( K \right)$ a cusp realizing the maximum in \eqref{eq:defMu1}, we also set
                \begin{equation}
                    \begin{array}[t]{lll}
                        \mu_{\mathfrak{a},2} \left( \tau \right) & = & \max\limits_{\substack{c \in \mathbb{P}^1 \left( K \right) \\ c \neq c_{\tau}}} \mu_{\mathfrak{a}} \left( \tau, \, c \right).
                    \end{array}
                \end{equation}
            \end{definition}

            \begin{proposition}
            \label{prop:mu1mu2NarrowClass}
                Let $\lambda \in K$ be totally positive. For any $\tau \in \mathbb{H}^n$, we have
                \begin{equation}
                    \begin{array}{lll}
                        \mu_{\lambda \mathfrak{a},j} \left( \tau \right) & = & \displaystyle \frac{1}{N \left( \lambda \right)} \mu_{\mathfrak{a},j} \left( \lambda \tau \right),
                    \end{array}
                \end{equation}
                for any $j \in \left \{ 1,2 \right \}$.
            \end{proposition}

            \begin{proof}
                Let us see the computation for $j = 1$. For any $\tau \in \mathbb{H}^n$, we have
                \begin{equation}
                    \begin{array}{llllll}
                        \multicolumn{3}{l}{\mu_{\lambda \mathfrak{a},1} \left( \tau \right)} & = & \max\limits_{\left[ \alpha : \beta \right] \in \mathbb{P}^1 \left( K \right)} \mu_{\lambda \mathfrak{a},1} \left( \tau, \left[ \alpha : \beta \right] \right) \\[1em]

                        & = & \multicolumn{3}{l}{\displaystyle \frac{1}{N \left( \lambda \right)} \max\limits_{\left[ \alpha : \beta \right] \in \mathbb{P}^1 \left( K \right)} \mu_{\mathfrak{a},1} \left( \lambda \tau, \left[ \alpha : \lambda^{-1} \beta \right] \right)} & \text{using proposition \ref{prop:muFunctionNarrowClass}} \\[1em]

                        & = & \multicolumn{3}{l}{\displaystyle \frac{1}{N \left( \lambda \right)} \max\limits_{\left[ \alpha : \beta \right] \in \mathbb{P}^1 \left( K \right)} \mu_{\mathfrak{a},1} \left( \lambda \tau, \left[ \alpha : \beta \right] \right)} \\[1em]

                        & = & \multicolumn{3}{l}{\displaystyle \frac{1}{N \left( \lambda \right)} \mu_{\mathfrak{a},1} \left( \lambda \tau \right).}
                    \end{array}
                \end{equation}
            \end{proof}

        \subsubsection{\texorpdfstring{The fundamental $\mathfrak{a}$-neighbourhoods of the cusps}{The fundamental a-neighbourhoods of the cusps}}

            As before, we consider a fractional ideal $\mathfrak{a}$ of $\mathcal{O}_K$.

            \begin{definition}
                For any $c \in \mathbb{P}^1 \left( K \right)$ and any $r>0$, we set
                \begin{equation}
                    \begin{array}{lll}
                        B_{\mathfrak{a}} \left( c, r \right) & = & \left \{ \tau \in \mathbb{H}^n \; \middle \vert \; \mu_{\mathfrak{a}} \left( \tau, c \right)^{-1/2} < r \right \},
                    \end{array}
                \end{equation}
                which can be thought of as the ball of $\mathfrak{a}$-radius $r$ centred at the cusp $c$.
            \end{definition}

            \begin{remark}
                When we take $c = \infty$, the fractional ideal $\mathfrak{a}$ becomes irrelevant in the definition above, as it was then in the definition of $\mu_{\mathfrak{a}} = \mu$. In this case, the notation $B_{\mathfrak{a}} \left( \infty, r \right)$ may be lightened to $B \left( \infty, r \right)$.
            \end{remark}

            \begin{proposition}
                For any $c \in \mathbb{P}^1 \left( K \right)$, the stabilizer $\widehat{\Gamma}_K \left( \mathfrak{a} \right)_c$ acts on $B_{\mathfrak{a}} \left( c, r \right)$.
            \end{proposition}

            \begin{proof}
                This is a direct consequence of proposition \ref{prop:functionMuActionGamma}.
            \end{proof}

            \begin{lemma}
            \label{lemma:aBallsCToInfinity}
                For any $r>0$ and any cusp $c \in \mathbb{P}^1 \left( K \right)$, the map
                \begin{equation}
                    \begin{array}{ccc}
                        \mathbb{H}^n & \longrightarrow & \mathbb{H}^n \\
                        \tau & \longmapsto & M \cdot \tau
                    \end{array},
                \end{equation}
                where $M$ and $\mathfrak{q}$ are as in proposition \ref{prop:distanceMatrixMCToInfitnity}, induces a diffeomorphism 
                \begin{equation}
                    \begin{array}{lll}
                        \widehat{\Gamma}_K \left( \mathfrak{a} \mathfrak{q^2} \right)_{\infty} \backslash B_{\mathfrak{a q^2}} \left( \infty, N \left( \mathfrak{q} \right) r \right) & \simeq & \widehat{\Gamma}_K \left( \mathfrak{a} \right)_c \backslash B_{\mathfrak{a}} \left( c, r \right)
                    \end{array}
                \end{equation}
                whose Jacobian equals $1$.
            \end{lemma}

            \begin{proof}
                This is a direct consequence of proposition \ref{prop:distanceMatrixMCToInfitnity}.
            \end{proof}

            \begin{proposition}
                For any $r>0$ and any cusp $c = \left[ \alpha : \beta \right] \in \mathbb{P}^1 \left( K \right)$, we have
                \begin{equation}
                    \begin{array}{lll}
                        \vol \left( \widehat{\Gamma}_K \left( \mathfrak{a} \right)_c \backslash B_{\mathfrak{a}} \left( c, r \right) \right) & = & \displaystyle N \left( \mathfrak{a} \right)^{-1} \sqrt{\Delta_K} \cdot \frac{2^{n-1}}{[ \mathcal{O}_K^{\times,+} : \mathcal{O}_K^{\times, 2}]} \cdot R_K \cdot r^2.
                    \end{array}
                \end{equation}
            \end{proposition}

            \begin{proof}
                The first step is to use lemma \ref{lemma:aBallsCToInfinity} to obtain
                \begin{equation}
                    \begin{array}{lll}
                        \vol \left( \widehat{\Gamma}_K \left( \mathfrak{a} \right)_c \backslash B_{\mathfrak{a}} \left( c, r \right) \right) & = & \vol \left( \widehat{\Gamma}_K \left( \mathfrak{a} \mathfrak{q^2} \right)_{\infty} \backslash B \left( \infty, N \left( \mathfrak{q} \right) r \right) \right),
                    \end{array}
                \end{equation}
                thereby reducing the computation to the cusp $\infty$. Note that we have
                \begin{equation}
                    \begin{array}{lllll}
                        \vol T & = & \vol \left( \mathfrak{b}^{-1} \backslash \mathbb{R}^n \right) & = & N \left( \mathfrak{b} \right)^{-1} \sqrt{\Delta_K},
                    \end{array}
                \end{equation}
                where, as in proposition \ref{prop:fundamentalDomainStabInfty}, we denote by $T$ a fundamental domain for the action of $\mathfrak{b}$ on $\mathbb{R}^n$ using the embeddings $\sigma_1, \ldots, \sigma_n$. A computation similar to the one made in the proof of \cite[Proposition 3.26]{dutour2025minkowskitypetheoremdistancescusps} then shows that, for any $t>0$, we have
                \begin{equation}
                    \begin{array}{lll}
                        \vol \left( \widehat{\Gamma}_K \left( \mathfrak{b} \right)_{\infty} \backslash B \left( \infty, t \right) \right) & = & \displaystyle N \left( \mathfrak{b} \right)^{-1} \sqrt{\Delta_K} \cdot \frac{2^{n-1}}{[ \mathcal{O}_K^{\times,+} : \mathcal{O}_K^{\times, 2}]} \cdot R_K \cdot t^2.
                    \end{array}
                \end{equation}
                Taking in particular $\mathfrak{b} = \mathfrak{a} \mathfrak{q}^2$ and $t = N \left( \mathfrak{q} \right) r$ yields
                \begin{equation}
                    \begin{array}{lll}
                        \vol \left( \widehat{\Gamma}_K \left( \mathfrak{a} \right)_c \backslash B_{\mathfrak{a}} \left( c, r \right) \right) & = & \vol \left( \widehat{\Gamma}_K \left( \mathfrak{a} \mathfrak{q^2} \right)_{\infty} \backslash B \left( \infty, N \left( \mathfrak{q} \right) r \right) \right) \\[2em]

                        & = & \displaystyle N \left( \mathfrak{a} \right)^{-1} \sqrt{\Delta_K} \cdot \frac{2^{n-1}}{[ \mathcal{O}_K^{\times,+} : \mathcal{O}_K^{\times, 2}]} \cdot R_K \cdot r^2,
                    \end{array}
                \end{equation}
                which completes the proof of the proposition.
            \end{proof}

            \begin{remark}
                In particular, we obtain the quadratic behaviour in $r$ of these volumes, as we have
                \begin{equation}
                    \begin{array}{lll}
                        \vol \left( \widehat{\Gamma}_K \left( \mathfrak{a} \right)_c \backslash B_{\mathfrak{a}} \left( c, r \right) \right) & = & r^2 \vol \left( \widehat{\Gamma}_K \left( \mathfrak{a} \right)_c \backslash B_{\mathfrak{a}} \left( c, 1 \right) \right).
                    \end{array}
                \end{equation}
            \end{remark}

            Let us conclude this section by defining the notion of $\mathfrak{a}$-spheres of influence and using them to obtain a fundamental domain for the action of $\widehat{\Gamma}_K \left( \mathfrak{a} \right)$ on $\mathbb{H}^n$.

            \begin{definition}
            \label{def:sphereInfluence}
                Let $c \in \mathbb{P}^1 \left( K \right)$. The \textit{$\mathfrak{a}$-sphere of influence} of $c$ is defined as
                \begin{equation}
                    \begin{array}{lll}
                        S_{\mathfrak{a}, c} & = & \left \{ \tau \in \mathbb{H}^n \; \middle \vert \; \mu_{\mathfrak{a}, 1} \left( \tau \right) \, = \, \mu_{\mathfrak{a}} \left( \tau, c \right) \right \},
                    \end{array}
                \end{equation}
                \textit{i.e.} as the set of points whose $\mathfrak{a}$-distance to $c$ is smaller than their $\mathfrak{a}$-distance to any other cusp.
            \end{definition}

            \begin{proposition}
                Let $\gamma \in \widehat{\Gamma}_K \left( \mathfrak{a} \right)$. We have
                \begin{equation}
                    \begin{array}{lll}
                        \gamma \cdot S_{\mathfrak{a},c}^{\circ} \cap S_{\mathfrak{a},c}^{\circ} \neq \emptyset & \iff & \gamma \in \widehat{\Gamma}_K \left( \mathfrak{a} \right)_c,
                    \end{array}
                \end{equation}
                where $S_{\mathfrak{a},c}^{\circ}$ denotes the interior of $S_{\mathfrak{a}, c}$.
            \end{proposition}

            \begin{proof}
                The argument is the same as the one used in \cite[Proposition I.2.3]{vanDerGeer:hilbert-modular-surfaces}.
            \end{proof}

            \begin{proposition}
            \label{prop:fundamentalDomain}
                Consider a set $c_1, \ldots, c_{h_K} \in \mathbb{P}^1 \left( K \right)$ of representatives of the cusps modulo the action of $\widehat{\Gamma}_K \left( \mathfrak{a} \right)$. The set
                \begin{equation}
                    \begin{array}{lll}
                        S & = & \bigsqcup\limits_{j=1}^{h_K} \widehat{\Gamma}_K \left( \mathfrak{a} \right)_{c_j} \backslash S_{\mathfrak{a}, c_j}
                    \end{array}
                \end{equation}
                is a fundamental domain for the action of $\widehat{\Gamma}_K \left( \mathfrak{a} \right)$ on $\mathbb{H}^n$.
            \end{proposition}

            \begin{proof}
                This is similar to the case $\mathfrak{a} = \mathcal{O}_K$, which is detailed in \cite[Section I.2]{vanDerGeer:hilbert-modular-surfaces}.
            \end{proof}

\section{A Minkowski-type theorem}

    In this final section, we will extend the results obtained in \cite{dutour2025minkowskitypetheoremdistancescusps} for totally real number fields whose class number equals $1$ to arbitrary totally real number fields. Some of the material presented in \cite{dutour2025minkowskitypetheoremdistancescusps} will be used again here, and we will go in further detail whenever appropriate.

    \subsection{\texorpdfstring{From points in $\mathbb{H}^n$ to rigid adelic spaces}{From points in Hn to rigid adelic spaces}}

        Recall that, in \cite[Section 4.1]{dutour2025minkowskitypetheoremdistancescusps}, we saw that the maps
        \begin{equation}
        \label{eq:correspondanceS2H}
            \begin{array}{ccccccc}
                \multicolumn{3}{c}{S_2^{++} \left( \mathbb{R} \right)} & \longrightarrow & \multicolumn{3}{c}{\mathbb{H} \times \mathbb{R}_+^{\ast}} \\[1em]

                S & = & \left( \begin{array}{cc} u & v \\[0.4em] v & w \end{array} \right) & \longmapsto & \displaystyle \left( \tau_S, \det S \right) & = & \displaystyle \left( \frac{v + i \sqrt{\det S}}{w}, \, \det S \right) \\[2em]

                \multicolumn{3}{c}{\displaystyle \frac{\sqrt{\lambda}}{y} \left( \begin{array}{cc} x^2+y^2 & x \\[0.4em] x & 1 \end{array} \right)} & \longmapsfrom & \multicolumn{3}{c}{\displaystyle \left( \tau \, = \, x + i y, \, \lambda \right)}
            \end{array}
        \end{equation}
        give a correspondence between $\mathbb{H} \times \mathbb{R}_+^{\ast}$ and $S_2^{++} \left( \mathbb{R} \right)$, which is compatible with the natural group actions of $PSL_2 \left( \mathbb{R} \right)$. In particular, to any $\tau \in \mathbb{H}$, we attach the matrix
        \begin{equation}
            \begin{array}{lllll}
                S \left( \tau \right) & = & \displaystyle \frac{1}{y} \left( \begin{array}{cc} x^2+y^2 & x \\[0.4em] x & 1 \end{array} \right) & = & \transpose{T \left( \tau \right)} T \left( \tau \right)
            \end{array}
        \end{equation}
        where we have written $\tau \, = \, x + iy$, and with
        \begin{equation}
            \begin{array}{lll}
                T \left( \tau \right) & = & \displaystyle \frac{1}{\sqrt{y}} \left( \begin{array}{cc} y & 0 \\[0.4em] x & 1 \end{array} \right).
            \end{array}
        \end{equation}
        After iteration, we obtain a correspondence
        \begin{equation}
            \begin{array}{lll}
                S_2^{++} \left( \mathbb{R} \right)^n & \simeq & \mathbb{H}^n \times \left( \mathbb{R}_+^{\ast} \right)^n,
            \end{array}
        \end{equation}
        which is compatible with the actions of $PSL_2 \left( \mathbb{R} \right)^n$, and in particular with the actions of any generalized Hilbert modular group $\widehat{\Gamma}_K \left( \mathfrak{a} \right)$, where $\mathfrak{a}$ is a fractional ideal of $\mathcal{O}_K$. Let $\mathfrak{a}$ be such an ideal, and write
        \begin{equation}
            \begin{array}{lll}
                \mathfrak{a} \mathcal{O}_{K_v} & = & x_v \mathcal{O}_{K_v}
            \end{array}
        \end{equation}
        for every place $v \in V_f \left( K \right)$.
        \begin{definition}
        \label{def:rigidAdelicSpaceEATau}
            For any point $\tau = \left( \tau_1, \ldots, \tau_n \right) \in \mathbb{H}^n$, we define the rigid adelic space $E_{\mathfrak{a}, \tau}$ as given by the matrices
            \begin{equation}
                \begin{array}{llll}
                    A_{\mathfrak{a}, v} & = & \left( \begin{array}{cc} 1 & 0 \\[0.4em] 0 & x_v \end{array} \right) & \text{if } v \text{ is finite}, \\[2em]
                    
                    A_{j,v} & = & T \left( \tau_j \right) & \text{if } v \text{ is infinite, associated with } \sigma_j ,
                \end{array}
            \end{equation}
            which together yield an element of $GL_2 \left( \mathbb{A}_K \right)$.
        \end{definition}

        \begin{proposition}
        \label{prop:HETauComputation}
            For any $\tau \in \mathbb{H}^n$, we have
            \begin{equation}
                \begin{array}{lll}
                    H \left( E_{\mathfrak{a}, \tau} \right) & = & N \left( \mathfrak{a} \right)^{-1/n}.
                \end{array}
            \end{equation}
        \end{proposition}

        \begin{proof}
            We have
            \begin{equation}
                \begin{array}{llll}
                    H \left( E_{\mathfrak{a}, \tau} \right) & = & \prod\limits_{v \in V_f \left( K \right)} \left \vert x_v \right \vert_v^{n_v/n} \cdot \prod\limits_{j=1}^n \left \vert \det T \left( \tau_j \right) \right \vert^{1/n} \\[2em]

                    & = & N \left( \mathfrak{a} \right)^{-1/n}
                \end{array}
            \end{equation}
            using the computation done in the proof of corollary \ref{cor:computeNormIdeal}.
        \end{proof}

        The following lemma, which is a corollary to proposition \ref{cor:computeNormIdeal}, will be useful later in the proof of proposition \ref{prop:heightAndMuFunction}, so we record it here.

        \begin{lemma}
        \label{lem:productNormsFinitePlaces}
            Let $\mathfrak{a}$ be an ideal of $\mathcal{O}_K$. For any $\alpha, \beta \in \mathcal{O}_K$, we have
            \begin{equation}
                \begin{array}{lll}
                    \prod\limits_{v \in V_f \left( K \right)} \left \Vert \left( \alpha, \beta \right) \right \Vert^{n_v/n}_{A_{\mathfrak{a}, v}, v} & = & N \left( \alpha \mathcal{O}_K + \beta \mathfrak{a} \right)^{-1/n}.
                \end{array}
            \end{equation}
        \end{lemma}

        \begin{proof}
            Consider $\alpha, \beta \in \mathcal{O}_K$. For a place $v \in V_f \left( K \right)$, we have
            \begin{equation}
                \begin{array}{lllll}
                    \left \Vert \left( \alpha, \beta \right) \right \Vert_{A_{\mathfrak{a}, v}, v} & = & \left \Vert A_{\mathfrak{a},v} \left( \begin{array}{c} \alpha \\ \beta \end{array} \right) \right\Vert_{v} & = & \max \left( \left \vert \alpha \right \vert_v, \left \vert \beta x_v \right \vert_v \right).
                \end{array}
            \end{equation}
            The proposition then becomes a direct consequence of corollary \ref{cor:computeNormIdeal}.
        \end{proof}

    \subsection{The relation between heights and distances to cusps}

        In this paragraph, we will see that the heights $H_{E_{\mathfrak{a}, \tau}}$ of the points in the rigid adelic space $E_{\mathfrak{a}, \tau}$ can be related to the $\mathfrak{a}$-distances, introduced in definition \ref{def:aDistance}, between $\tau$ and the cusps, whose set is identified with $\mathbb{P}^1 \left( K \right)$. The differences with \cite{dutour2025minkowskitypetheoremdistancescusps} are that the class number assumption $h_K = 1$ is not made here, and that we consider arbitrary fractional ideals.

        \begin{definition}
            We define the involution $\iota$ of $\mathbb{P}^1 \left( K \right)$ by
            \begin{equation}
                \begin{array}{ccccc}
                    \iota & : & \mathbb{P}^1 \left( K \right) & \longrightarrow & \mathbb{P}^1 \left( K \right) \\[0.5em]
                    && \left[ x,y \right] & \longmapsto & \left[ y : -x \right]
                \end{array}.
            \end{equation}
        \end{definition}
    
        \begin{proposition}
        \label{prop:heightAndMuFunction}
            For any fractional ideal $\mathfrak{a}$ and any $\tau \in \mathbb{H}^n$, we have
            \begin{equation}
            \label{eq:heightAndMuFunction}
                \begin{array}{ccc}
                    H_{E_{\mathfrak{a}, \tau}} \circ i & = & N \left( \mathfrak{a} \right)^{-1/n} \mu_{\mathfrak{a}} \left( \tau, \, \cdot \right)^{-1/2n},
                \end{array}
            \end{equation}
            where $\mu_{\mathfrak{a}}$ is the function introduced in definition \ref{def:aDistance}.
        \end{proposition}

        \begin{proof}
            Consider a fractional ideal $\mathfrak{a}$, as well as points $\tau \in \mathbb{H}^n$ and $\left[ \alpha : \beta \right] \in \mathbb{P}^1 \left( K \right)$. We have
            \begin{equation}
                \begin{array}{llllll}
                    \multicolumn{3}{l}{H_{E_{\mathfrak{a}, \tau}} \left( \left[ \alpha : \beta \right] \right)} & = & \multicolumn{2}{l}{\prod\limits_{v \in V_f \left( K \right)} \left \Vert \left( \alpha, \beta \right) \right \Vert^{n_v/n}_{A_{\mathfrak{a}, v}, v} \cdot \prod\limits_{j=1}^n \left \Vert \left( \alpha, \beta \right) \right \Vert^{1/n}_{\tau,j}} \\[2em]

                    & = & \multicolumn{3}{l}{N \left( \alpha \mathcal{O}_K + \beta \mathfrak{a} \right)^{-1/n} \cdot \prod\limits_{j=1}^n \left \Vert \left( \alpha, \beta \right) \right \Vert^{1/n}_{\tau,j}} & \text{by lemma \ref{lem:productNormsFinitePlaces}} \\[2em]

                    & = & \multicolumn{4}{l}{N \left( \alpha \mathcal{O}_K + \beta \mathfrak{a} \right)^{-1/n} \cdot \prod\limits_{j=1}^n y_j^{-1/2n} \left \vert \sigma_j \left( \alpha \right) \tau_j + \sigma_j \left( \beta \right) \right \vert^{1/n}} \\[2em]

                    & = & \multicolumn{3}{l}{N \left( \mathfrak{a} \right)^{-1/n} \mu_{\mathfrak{a}} \left( \tau, \iota \left( \left[ \alpha : \beta \right] \right) \right)^{-1/2n},}
                \end{array}
            \end{equation}
            thereby concluding the proof, since $\iota$ is an involution.
        \end{proof}

        \begin{remark}
            Note that taking $\mathfrak{a} = \mathcal{O}_K$, which satisfies $N \left( \mathfrak{a} \right) = 1$ and $h_K = 1$ yields the same result as \cite[Proposition 4.10]{dutour2025minkowskitypetheoremdistancescusps}.
        \end{remark}

        As in \cite[Proposition 4.12]{dutour2025minkowskitypetheoremdistancescusps}, we will now use proposition \ref{prop:heightAndMuFunction} to relate the Roy--Thunder minima of $E_{\mathfrak{a}, \tau}$ to the $\mathfrak{a}$-distances between $\tau$ and its two closest cusps.

        \begin{proposition}
        \label{prop:RoyThunderMinimaMu}
            For any fractional ideal $\mathfrak{a}$ and any $\tau \in \mathbb{H}^n$, we have
            \begin{equation}
            \label{eq:Lambda1Lambda2Mu}
                \begin{array}{lll}
                    \Lambda_1 \left( E_{\mathfrak{a}, \tau} \right) & = & N \left( \mathfrak{a} \right)^{-1/n} \mu_{\mathfrak{a}, 1} \left( \tau \right)^{-1/2n}, \\[0.5em]

                    \Lambda_2 \left( E_{\mathfrak{a}, \tau} \right) & = & N \left( \mathfrak{a} \right)^{-1/n} \mu_{\mathfrak{a}, 2} \left( \tau \right)^{-1/2n}.
                \end{array}
            \end{equation}
        \end{proposition}

        \begin{proof}
            The proof is entirely similar to that of \cite[Proposition 4.12]{dutour2025minkowskitypetheoremdistancescusps}.
        \end{proof}

    \subsection{A Minkowski-type theorem}

        Let us now use proposition \ref{prop:RoyThunderMinimaMu} to obtain a version of Minkowski's second theorem on the function $\mu_{\mathfrak{a}}$.

        \begin{theorem}
        \label{thm:minkowski}
            For any fractional ideal $\mathfrak{a}$ and any $\tau \in \mathbb{H}^n$, we have
            \begin{equation}
                \begin{array}{lllll}
                    \displaystyle \frac{1}{c_{II}^{\Lambda} \left( 2, K \right)^{4n}} \cdot \frac{1}{N \left( \mathfrak{a} \right)^2} & \leqslant & \mu_{\mathfrak{a}, 1} \left( \tau \right) \mu_{\mathfrak{a}, 2} \left( \tau \right) & \leqslant & \displaystyle \frac{1}{N \left( \mathfrak{a} \right)^2},
                \end{array}
            \end{equation}
            where $c_{II}^{\Lambda} \left( 2, K \right)$ is the constant introduced in \cite[Definition 2.31]{dutour2025minkowskitypetheoremdistancescusps}.
        \end{theorem}

        \begin{proof}
            This is a direct consequence of propositions \ref{prop:HETauComputation}, \ref{prop:RoyThunderMinimaMu}, and of \cite[Proposition 2.33]{dutour2025minkowskitypetheoremdistancescusps}.
        \end{proof}

        Let us see two consequences of theorem \ref{thm:minkowski}, which generalize \cite[Corollary 4.15, 4.17]{dutour2025minkowskitypetheoremdistancescusps}.

        \begin{corollary}
        \label{cor:separationCusps}
            Let $c, c' \in \mathbb{P}^1 \left( K \right)$ be two cusps, and $\mathfrak{a}$ be a fractional ideal. For any $\tau \in \mathbb{H}^n$, we have
            \begin{equation}
                \begin{array}{lll}
                    \displaystyle \mu_{\mathfrak{a}} \left( \tau, c \right) \; \geqslant \; \frac{1}{N \left( \mathfrak{a} \right)} \text{ and } \mu_{\mathfrak{a}} \left( \tau, c' \right) \; \geqslant \; \frac{1}{N \left( \mathfrak{a} \right)} & \implies & c \; = \; c'.
                \end{array}
            \end{equation}
            In other words, a point in $\mathbb{H}^n$ cannot be at an $\mathfrak{a}$-distance less than $N \left( \mathfrak{a} \right)^{1/2}$ from two different cusps. As a consequence, we have
            \begin{equation}
            \label{eq:inclusionBInfinity1IntoSInfinity}
                \begin{array}{lll}
                    B_{\mathfrak{a}} \left( c, N \left( \mathfrak{a} \right)^{1/2} \right) & \subseteq & S_{\mathfrak{a}, c} ,
                \end{array}
            \end{equation}
            where we recall that $S_{\mathfrak{a}, c}$ is the $\mathfrak{a}$-sphere of influence of $c$, and $B_{\mathfrak{a}} \left( c, r \right)$ is the ball centred at $c$ of radius $r > 0$.
        \end{corollary}

        \begin{proof}
            By contraposition, assuming we have $c \neq c'$, the inequality
            \begin{equation}
                \begin{array}{lll}
                    \mu_{\mathfrak{a}, 1} \left( \tau \right) \mu_{\mathfrak{a}, 2} \left( \tau \right) & \leqslant & \displaystyle \frac{1}{N \left( \mathfrak{a} \right)^2}
                \end{array}
            \end{equation}
            and the fact that, by definition, we have $\mu_{\mathfrak{a}, 1} \left( \tau \right) \, \geqslant \, \mu_{\mathfrak{a}, 2} \left( \tau \right)$, give
            \begin{equation}
                \begin{array}{lll}
                    \mu_{\mathfrak{a}, 2} \left( \tau \right) & \leqslant & \displaystyle \frac{1}{N \left( \mathfrak{a} \right)}.
                \end{array}
            \end{equation}
            Hence, since we have assumed $c \neq c'$, we must have 
            \begin{equation}
                \begin{array}{lll}
                   \displaystyle  \mu_{\mathfrak{a}} \left( \tau, c \right) \, < \, \frac{1}{N \left( \mathfrak{a} \right)} & \text{or} & \displaystyle \mu_{\mathfrak{a}} \left( \tau, c' \right) \, < \, \frac{1}{N \left( \mathfrak{a} \right)},
                \end{array}
            \end{equation}
            thus completing the proof.
        \end{proof}

        \begin{remark}
        \label{rmk:comparisonVDG}
            Compared to \cite[Lemma I.2.1]{vanDerGeer:hilbert-modular-surfaces}, note that corollary \ref{cor:separationCusps} is more general, as it takes into account generalized Hilbert modular groups $\widehat{\Gamma}_K \left( \mathfrak{a} \right)$, effective, and uniform, in the sense that the dependance on $K$ only comes from the choice of the fractional ideal $\mathfrak{a}$.
        \end{remark}

        \begin{corollary}
        \label{cor:lowerBoundMu1}
            Consider a fractional ideal $\mathfrak{a}$. For any $\tau \in \mathbb{H}^n$, we have
            \begin{equation}
            \label{eq:lowerBoundMu1}
                \begin{array}{lll}
                    \displaystyle \mu_{\mathfrak{a}, 1} \left( \tau \right) & \geqslant & \displaystyle \frac{1}{c_{II}^{\Lambda} \left( 2, K \right)^{2n}} \cdot \frac{1}{N \left( \mathfrak{a} \right)}.
                \end{array}
            \end{equation}
            In other words, a point in $\mathbb{H}^n$ is always at an $\mathfrak{a}$-distance at most $c_{II}^{\Lambda} \left( 2, K \right)^n N \left( \mathfrak{a} \right)^{1/2}$ from a cusp. As a consequence, for any cusp $c \in \mathbb{P}^1 \left( K \right)$, we have
            \begin{equation}
            \label{eq:inclusionSInfinityIntoBInfinityC2K}
                \begin{array}{lll}
                    S_{\mathfrak{a}, c} & \subseteq & B \left( c, c_{II}^{\Lambda} \left( 2, K \right)^n N \left( \mathfrak{a} \right)^{1/2} \right).
                \end{array}
            \end{equation}
        \end{corollary}

        \begin{proof}
            The inequality
            \begin{equation}
                \begin{array}{lllll}
                    \displaystyle \frac{1}{c_{II}^{\Lambda} \left( 2, K \right)^{4n}} \cdot \frac{1}{N \left( \mathfrak{a} \right)^2}& \leqslant & \mu_{\mathfrak{a}, 1} \left( \tau \right) \mu_{\mathfrak{a}, 2} \left( \tau \right) & \leqslant & \mu_{\mathfrak{a}, 1} \left( \tau \right)^2
                \end{array}
            \end{equation}
            directly gives the result.
        \end{proof}

        \begin{remark}
            Compared to \cite[Lemma I.2.1]{vanDerGeer:hilbert-modular-surfaces}, note that corollary \ref{cor:lowerBoundMu1} is more general, effective, and uniform, for the reasons given in remark \ref{rmk:comparisonVDG}, and that it is also optimal, by definition of the constant $c_{II}^{\Lambda} \left( 2, K \right)$. Using this for $\mathfrak{a} \, = \, \mathcal{O}_K$ and the computation made by van der Geer in the proof of \cite[Lemma I.2.2]{vanDerGeer:hilbert-modular-surfaces}, we obtain the estimate
            \begin{equation}
                \begin{array}{lll}
                    c_{II}^{\Lambda} \left( 2, K \right) & \leqslant & \sqrt{2} \Delta_K^{1/2n},
                \end{array}
            \end{equation}
            which, as was already noticed in \cite{dutour2025minkowskitypetheoremdistancescusps}, is the same estimate as the one obtained by Gaudron and Rémond in \cite[Proposition 5.1]{gaudron-remond:corps-siegel}.
        \end{remark}

        Even though the following statement is a particular case of corollary \ref{cor:lowerBoundMu1}, let us explicitly state it, as it deals with the codifferent ideal $\mathfrak{d}_K^{-1}$, which is the one considered in \cite{frey-lefourn-lorenzo:height-estimates, habegger-pazuki:bad-reduction}.

        \begin{corollary}
            For any $\tau \in \mathbb{H}^n$, we have
            \begin{equation}
                \begin{array}{lll}
                    \mu_{\mathfrak{d}_K^{-1}, 1} \left( \tau \right) & \geqslant & 2^{-n},
                \end{array}
            \end{equation}
            where $\mathfrak{d}_K$ is the different ideal (see definition \ref{def:differentIdeal}).
        \end{corollary}

        \begin{proof}
            Taking $\mathfrak{a} = \mathfrak{d}_K^{-1}$, and using the estimate
            \begin{equation}
                \begin{array}{lll}
                    c_{II}^{\Lambda} \left( 2, K \right) & \leqslant & \sqrt{2} \Delta_K^{1/2n}
                \end{array}
            \end{equation}
            gives this result as a direct consequence of corollary \ref{cor:lowerBoundMu1}.
        \end{proof}

        \begin{proposition}
            For any cusp $c \in \mathbb{P}^1 \left( K \right)$, the boundary $\partial S_{\mathfrak{a}, c}$ of the $\mathfrak{a}$-sphere of influence $S_{\mathfrak{a}, c}$ satisfies
            \begin{equation}
                \begin{array}{lll}
                    \partial S_{\mathfrak{a}, c} & \subset & B \left( c, c_{II}^{\Lambda} \left( 2, K \right)^n N \left( \mathfrak{a} \right)^{1/2}\right) \setminus B \left( c, N \left( \mathfrak{a} \right)^{1/2} \right).
                \end{array}
            \end{equation}
        \end{proposition}

          \begin{proof}
                Recall that, in definition \ref{def:sphereInfluence}, the $\mathfrak{a}$-sphere of influence $S_{\mathfrak{a}, c}$ was defined as
                \begin{equation}
                    \begin{array}{lll}
                        S_{\mathfrak{a}, c} & = & \left \{ \tau \in \mathbb{H}^n \; \middle \vert \; \mu_{\mathfrak{a}} \left( \tau, c \right) \, = \, \mu_{\mathfrak{a}, 1} \left( \tau \right) \right \}.
                    \end{array}
                \end{equation}
                The boundary $\partial S_{\mathfrak{a}, c}$ is then the set of point in $S_{\mathfrak{a},c}$ which are equidistant from $c$ and another cusp, which means that we have
                \begin{equation}
                    \begin{array}{lll}
                        \partial S_{\mathfrak{a}, c} & = & \left \{ \tau \in \mathbb{H}^n \; \middle \vert \; \mu_{\mathfrak{a}} \left( \tau, c \right) \, = \, \mu_{\mathfrak{a}, 1} \left( \tau \right) \, = \, \mu_{\mathfrak{a}, 2} \left( \tau \right) \right \}
                    \end{array}
                \end{equation}
                By theorem \ref{thm:minkowski}, a point $\tau \in \partial S_{\mathfrak{a}, c}$ must satisfy
                \begin{equation}
                \label{eq:ineqMuBoundary}
                    \begin{array}{lllll}
                        \displaystyle \frac{1}{c_{II}^{\Lambda} \left( 2, K \right)^{4n}} \cdot \frac{1}{N \left( \mathfrak{a} \right)^2} & \leqslant & \mu_{\mathfrak{a}} \left( \tau, c \right)^2 & \leqslant & \displaystyle \frac{1}{N \left( \mathfrak{a} \right)^2},
                    \end{array}
                \end{equation}
                which yields the result, after applying an exponent $-1/4$ to each part of \eqref{eq:ineqMuBoundary}.
          \end{proof}

    \subsection{An interesting class of integrals}

        Simlarly to \cite[Section 4.5]{dutour2025minkowskitypetheoremdistancescusps}, we will conclude this paper by obtaining estimates on a class of integrals of the form
        \begin{equation}
        \label{eq:interestingIntegrals}
            \begin{array}{lll}
                \displaystyle \frac{1}{\vol \left( \widehat{\Gamma}_K \left( \mathfrak{a} \right) \backslash \mathbb{H}^n \right)} \int_{\widehat{\Gamma}_K \left( \mathfrak{a} \right) \backslash \mathbb{H}^n} \mu_{\mathfrak{a},1} \left( \tau \right)^t \; \mathrm{d}m\left( \tau \right),
            \end{array}
        \end{equation}
        with $t$ to be adjusted later, and $\mathfrak{a}$ a fractional ideal of $\mathcal{O}_K$. 

        \subsubsection{The partial volume function}

            Let us consider, as we did in proposition \ref{prop:fundamentalDomain}, representatives $c_1, \ldots, c_{h_K}$ of the cusps modulo the action of $\widehat{\Gamma}_K \left( \mathfrak{a} \right)$, so that
            \begin{equation}
                \begin{array}{lll}
                    S & = & \bigsqcup\limits_{j=1}^{h_K} \widehat{\Gamma}_K \left( \mathfrak{a} \right)_{c_j} \backslash S_{\mathfrak{a}, c_j}
                \end{array}
            \end{equation}
            is a fundamental domain for the action of $\widehat{\Gamma}_K \left( \mathfrak{a} \right)$ on $\mathbb{H}^n$.

            \begin{definition}
                The \textit{partial volume function} is defined as
                \begin{equation}
                    \begin{array}{ccccc}
                        g & : & \mathbb{R}_+^{\ast} & \longrightarrow & \mathbb{R}_+^{\ast} \\[0.5em]
                        && x & \longmapsto & \sum\limits_{j=1}^{h_K} \vol \left( \widehat{\Gamma}_K \left( \mathfrak{a} \right)_{c_j} \backslash \left( S_{\mathfrak{a},c_j} \cap B_{\mathfrak{a}} \left( c_j, x \right) \right) \right)
                    \end{array}
                \end{equation}
            \end{definition}

            \begin{proposition}
                The partial volume function $g$ is increasing, and satisfies
                \begin{equation}
                    \begin{array}{llll}
                        g \left( x \right) & = & \displaystyle N \left( \mathfrak{a} \right)^{-1} \sqrt{\Delta_K} \cdot \frac{2^{n-1}}{[ \mathcal{O}_K^{\times,+} : \mathcal{O}_K^{\times, 2}]} \cdot R_K \cdot x^2 & \text{if } x \leqslant N \left( \mathfrak{a} \right)^{1/2}
                    \end{array}
                \end{equation}
                as well as
                \begin{equation}
                    \begin{array}{llll}
                        g \left( x \right) & = & \displaystyle \vol \left( \widehat{\Gamma}_K \left( \mathfrak{a} \right) \backslash \mathbb{H}^n \right) & \text{if } x \geqslant c_{II}^{\Lambda} \left( 2, K \right)^n N \left( \mathfrak{a} \right)^{1/2}.
                    \end{array}
                \end{equation}
            \end{proposition}

            \begin{proof}
                These formulas are direct consequences of corollaries \ref{cor:separationCusps} and \ref{cor:lowerBoundMu1}.
            \end{proof}

        \subsubsection{Estimating the integrals}

            In this last paragraph, we will state the results leading to the lower- and upper- bounds of the integrals \eqref{eq:interestingIntegrals}. The proofs will be omitted, though appropriate references will be made, since the arguments are entirely similar to the ones used in \cite[Section 4.5]{dutour2025minkowskitypetheoremdistancescusps}.
            \begin{proposition}
                For any $0 \leqslant t < 1$, we have
                \begin{equation}
                    \begin{array}{lll}
                        \multicolumn{3}{l}{\displaystyle \int_{\widehat{\Gamma}_K \left( \mathfrak{a} \right) \backslash \mathbb{H}^n} \mu_{\mathfrak{a}, 1} \left( \tau \right)^t \mathrm{d}m \left( \tau \right)} \\[2em]

                        \qquad \quad & = & \displaystyle \frac{\vol \left( \widehat{\Gamma}_K \left( \mathfrak{a} \right) \backslash \mathbb{H}^n \right)}{c_{II}^{\Lambda} \left( 2, K \right)^{2nt} N \left( \mathfrak{a} \right)^t} + 2t \int_0^{c_{II}^{\Lambda} \left( 2, K \right)^n N \left( \mathfrak{a} \right)^{1/2}} g \left( x \right) x^{-2t-1} \mathrm{d}x.
                    \end{array}
                \end{equation}
            \end{proposition}
            \begin{proof}
                This is similar to \cite[Proposition 4.24]{dutour2025minkowskitypetheoremdistancescusps}.
            \end{proof}
            \begin{proposition}
                For any $0 \leqslant t < 1$, we have
                \begin{equation}
                    \begin{array}{lll}
                        \multicolumn{3}{l}{\displaystyle \left( \frac{1}{1-t} - \frac{1}{c_{II}^{\Lambda} \left( 2, K \right)^{2nt}} \right) \cdot \frac{g \left( N \left( \mathfrak{a} \right)^{1/2} \right)}{N \left( \mathfrak{a} \right)^t}} \\[2em] 
                        
                        \quad & \leqslant & \displaystyle 2t \int_0^{c_{II}^{\Lambda} \left( 2, K \right)^n N \left( \mathfrak{a} \right)^{1/2}} g \left( x \right) x^{-2t-1} \mathrm{d}x \\[2em]

                        & \leqslant & \displaystyle \frac{t}{1-t} \cdot \frac{g \left( N \left( \mathfrak{a} \right)^{1/2} \right)}{N \left( \mathfrak{a} \right)^t} + \frac{\vol \left( \widehat{\Gamma}_K \left( \mathfrak{a} \right) \backslash \mathbb{H}^n \right)}{N \left( \mathfrak{a} \right)^t} \left( 1 - \frac{1}{c_{II}^{\Lambda} \left( 2, K \right)^{2nt}} \right). 
                    \end{array}
                \end{equation}
            \end{proposition}
            \begin{proof}
                This is similar to \cite[Proposition 4.25]{dutour2025minkowskitypetheoremdistancescusps}.
            \end{proof}
            \begin{theorem}
            \label{thm:estimatesIntegral}
                For any $0 \leqslant t < 1$, we have
                \begin{equation}
                    \begin{array}{lllll}
                        \multicolumn{5}{l}{\displaystyle \frac{1}{N \left( \mathfrak{a} \right)^t} \left( \frac{1}{c_{II}^{\Lambda} \left( 2, K \right)^{2nt}} \left( 1 - \frac{1}{c_{II}^{\Lambda} \left( 2, K \right)^{2n}} \right) + \frac{1}{1-t} \cdot \frac{1}{c_{II}^{\Lambda} \left( 2, K \right)^{2n}} \right)} \\[2em]

                        & \leqslant & \displaystyle \frac{1}{\vol \left( \widehat{\Gamma}_K \left( \mathfrak{a} \right) \backslash \mathbb{H}^n \right)} \int_{\widehat{\Gamma}_K \left( \mathfrak{a} \right) \backslash \mathbb{H}^n} \mu_{\mathfrak{a},1} \left( \tau \right)^t \; \mathrm{d}m\left( \tau \right) & \leqslant & \displaystyle \frac{1}{1-t} \cdot \frac{1}{N \left( \mathfrak{a} \right)^t}.
                    \end{array}
                \end{equation}
            \end{theorem}

            \begin{proof}
                The proof is similar to the one of \cite[Theorem 4.26]{dutour2025minkowskitypetheoremdistancescusps}.
            \end{proof}

            \begin{remark}
                Using propositions \ref{prop:multByLambdaMap} and \ref{prop:mu1mu2NarrowClass}, we note that the inequalities in theorem \ref{thm:estimatesIntegral} do not provide more information when $\mathfrak{a}$ is replaced with another representative in its narrow class.
            \end{remark}


\bibliographystyle{amsplain}

\bibliography{bibliography}

\end{document}